\documentclass[11pt]{article}
\usepackage[latin1]{inputenc}
\usepackage{amsmath,amsthm,amssymb}
\usepackage{amsfonts}
\usepackage{amsmath,amsthm,amssymb,amscd}
\usepackage{latexsym}
\usepackage{color}
\usepackage{graphicx}
\usepackage{mathrsfs}
\usepackage{cite}

\usepackage{color,enumitem,graphicx}
\usepackage[colorlinks=true,urlcolor=black,
citecolor=black,linkcolor=black,linktocpage,pdfpagelabels,
bookmarksnumbered,bookmarksopen]{hyperref}

\makeatletter \@addtoreset{equation}{section} \makeatother

\newtheorem{theorem}{Theorem}[section]

\newtheorem{proposition}{Proposition}[section]
\newtheorem{lemma}{Lemma}[section]
\newtheorem{remark}{Remark}[section]

\allowdisplaybreaks

\begin{document}
\title{A weighted Adams' inequality in $\mathbb{R}^4$ and applications}

\author{Wenjing Chen\footnote{Corresponding author.}\ \footnote{E-mail address:\, {\tt wjchen@swu.edu.cn} (W. Chen), {\tt sqzhangmaths@163.com} (S. Zhang).}\  \ and Shiqi Zhang\\
\footnotesize  School of Mathematics and Statistics, Southwest University,
Chongqing, 400715, P.R. China}

\date{ }
\maketitle

\begin{abstract}
{In this paper, we establish a weighted Adams' inequality in some appropriate weighted Sobolev space in $\mathbb{R}^4$. Then we give an improvement inequality by proving the concentration-compactness result. In the last part, we consider an application to elliptic equation involving new exponential growth at infinity in $\mathbb{R}^4$.}


\end{abstract}

\section{Introduction and statement of main results}

Let $\Omega\subset\mathbb{R}^N$, $N\geq2$ be a bounded domain. By Sobolev embedding theorem, we know that
$W^{k,p}_{0}(\Omega)\hookrightarrow L^q(\Omega)$, for $1\leq q\leq\frac{Np}{N-kp}$, $kp<N$ and $W^{k,p}_{0}(\Omega)\hookrightarrow L^q(\Omega)$, for $1\leq q<\infty$, $kp=N$. However, it is known that $W^{k,p}_{0}(\Omega)\hookrightarrow L^{\infty}(\Omega)$ does not hold. In this case, it was provided independently by Yudovich \cite{VIY}, Peetre \cite{JP}, Pohazaev \cite{SIP}, Trudinger \cite{NST} that $W^{1,N}_{0}(\Omega)\hookrightarrow L_{\phi}(\Omega)$ where $L_{\phi}(\Omega)$ is the Orlicz space with the N-function $\phi(t)=e^{|t|^{\frac{N}{N-1}}-1}$. To be precise, Moser \cite{JM} proved the following inequality:
\begin{equation*}
  \sup_{u\in W^{1,N}_{0}(\Omega),|\nabla u|_{L^N(\Omega)}\leq1} \int_{\Omega}e^{\alpha u^{\frac{N}{N-1}}}\mathrm{d}x<+\infty
  \quad \Longleftrightarrow \quad \alpha\leq\alpha_N=N\omega^{\frac{1}{N-1}}_{N-1},
\end{equation*}
where $\omega_{N-1}$ denotes the area of the unit sphere in $\mathbb{R}^N$.

This result has led to many related results. For instance, the attainability of the supremum has been proved by Carleson and Chang \cite{CC} when the domain is a unit ball in $\mathbb{R}^N$, by Flucher \cite{MF} on an open and bounded domain of $\mathbb{R}^2$ with smooth boundary and by Lin \cite{LKC} on bounded smooth domain in any dimension. For the unbounded domain, several extensions of Trudinger-Moser type inequalities have been established, see \cite{Cao,BR,BRL,LZ}. In particular, Ruf \cite{BR} proved the following sharp inequality:
\begin{equation*}
  \sup\limits_{u\in H^1(\mathbb{R}^2),\int_{\mathbb{R}^2}\left(|u|^2+|\nabla u|^2\right)\mathrm{d}x\leq1}
  \int_{\mathbb{R}^2}\left(e^{\alpha u^2}-1\right)\mathrm{d}x<+\infty\quad \Longleftrightarrow \quad\alpha\leq4\pi.
\end{equation*}
Li and Ruf \cite{BRL} generalized the inequality to the N-dimensional($N\geq2$) case as follows:
\begin{equation*}
  \sup\limits_{u\in W^{1,N}(\mathbb{R}^N),\int_{\mathbb{R}^N}\left(|u|^N+|\nabla u|^N\right)\mathrm{d}x\leq1}
  \int_{\mathbb{R}^N}\Phi_N\left(\alpha|u|^{\frac{N}{N-1}}\right)\mathrm{d}x<+\infty\quad \Longleftrightarrow \quad\alpha\leq\alpha_N,
\end{equation*}
where $\Phi_N(t)=e^t-\sum^{N-2}\limits_{j=0}\frac{t^j}{j!}$.
Trudinger-Moser inequalities were also given on manifolds \cite{F,L,Y}, on Heisenberg groups \cite{CL,LL}, without boundary conditions \cite{C,C2} and references therein.

Adams \cite{A} studied the similar result to higher order Sobolev spaces. Indeed, Adams proved that
\begin{equation}\label{boundbi}
  \sup\limits_{u\in W^{2,2}_{0}(\Omega),\int_{\Omega}|\triangle u|^2\mathrm{d}x\leq1}\int_{\Omega}e^{32\pi^2u^2}\mathrm{d}x<\infty.
\end{equation}
When $\Omega$ is replaced by the whole space $\mathbb{R}^4$, the integrals in (\ref{boundbi}) are infinite. But Ruf and Sani \cite{RS} established the corresponding Adams type inequality as follows:
\begin{equation*}
  \sup\limits_{u\in W^{2,2}_{rad}(\mathbb{R}^4),\|u\|_{W^{2,2}}\leq1}
  \int_{\mathbb{R}^4}\left(e^{32\pi^2u^2}-1\right)\mathrm{d}x<+\infty,
\end{equation*}
where $\|u\|^2_{W^{2,2}}=\int_{\mathbb{R}^4}\left(|\triangle u|^2+|\nabla u|^2+u^2\right)\mathrm{d}x$. More related results can see \cite{LL2,CLZ} and references therein.

Recently, many Trudinger-Moser type and Adams' type inequalities on weighted Sobolev spaces have been studied. Calanchi and Ruf \cite{CR} first studied the case of weighted Sobolev norms with weights of logarithmic type. In fact, let $B=B_1(0)$ be the unit open ball in $\mathbb{R}^2$, $\omega_0(x)=\left(\log\frac{e}{|x|}\right)^{\beta}$, $0\leq\beta\leq1$. Denote by $H^{1}_{0,rad}(B,\omega_0)$ the standard weighted radial Sobolev space defined as the space of all radial functions of the completion of $C^{\infty}_{0}$ with respect to the norm
\begin{equation*}
  \|u\|^2_{\omega_0}=\int^{}_{B}|\nabla u|^2\omega_0(x)\mathrm{d}x.
\end{equation*}
Then, for $0\leq\beta<1$,
\begin{equation*}
  \sup_{u\in H^{1}_{0,rad}(B,\omega_0),\|u\|_{\omega_0}\leq1}\int^{}_{B}e^{\alpha|u|^{\frac{2}{1-\beta}}}\mathrm{d}x<+\infty
  \quad \Longleftrightarrow \quad \alpha\leq\alpha_{2,\beta}=2\left[2\pi(1-\beta)\right]^{\frac{1}{1-\beta}},
\end{equation*}
for $\beta=1$,
\begin{equation*}
  \sup_{u\in H^{1}_{0,rad}(B,\omega_{0}),\|u\|_{\omega_{0}}\leq1}\int^{}_{B}e^{\alpha e^{2\pi u^2}}\mathrm{d}x<+\infty
  \quad \Longleftrightarrow \quad \alpha\leq2.
\end{equation*}

Moreover, Calanchi and Ruf \cite{CR2} generalized the similar results to the N-dimensional case.
Here, we have to mention the work where Nguyen \cite{NVH} provided a simple proof which based on a suitable change of function and obtained the existence of maximizers for this inequality when $\beta$ is small enough. The attainability of the supremum of this inequality has been given by Roy \cite{PR,PR2}. For the related results on unbounded domains, we can cite \cite{SA,SAJ,SR,SRn}. In particular, Aouaoui and Jlel \cite{SR} extended the work of Calanchi and Ruf to the whole space $\mathbb{R}^2$, considering the following weight:
\begin{equation}\label{omega1}
  \omega_1(x)=\begin{cases}
  \big(\log(\frac{e}{\mid x\mid})\big)^{\beta},&\ \text{if $|x|<1$},\\
  \chi_0(|x|),&\ \text{if $|x|\geq1$},
  \end{cases}
\end{equation}
where $0<\beta\leq1$ and $\chi_0:[1,+\infty)\rightarrow(0,+\infty)$ is a continuous function such that $\chi_0(1)=1$ and $\inf\limits_{t\geq1}\chi_0(t)>0$. Denote by $F_\beta$ as the space of all radial functions of the completion of $C^{\infty}_0(\mathbb{R}^2)$ with respect to the norm
\begin{equation*}
  \|u\|^2_{\omega_1}=\int_{\mathbb{R}^2}\left(|\nabla u|^2\omega_1(x)+|u|^2\right)\mathrm{d}x.
\end{equation*}
More precisely,

\noindent\textbf{Theorem A}
(i) Let $0<\beta<1$ and $\omega_1$ be defined by (\ref{omega1}). For all $u\in F_\beta$, we have
\begin{equation*}
  \int_{\mathbb{R}^2}\left(e^{|u|^{\frac{2}{1-\beta}}}-1\right)\mathrm{d}x<+\infty.\notag
\end{equation*}
Moreover,
\begin{equation*}
  \sup_{u\in F_\beta,\|u\|_{\omega_1}\leq1} \int_{\mathbb{R}^2}\left(e^{\alpha|u|^{\frac{2}{1-\beta}}}-1\right)\mathrm{d}x<+\infty,
\end{equation*}
if and only if
$
  \alpha<\alpha_{2,\beta}.
$

(ii) Let $\beta=1$ and $\omega_1$ be defined by (\ref{omega1}). For all $u\in F_1$, we have
\begin{equation*}
  \int_{\mathbb{R}^2}\left(e^{\alpha\left(e^{u^2}-1\right)}-1\right)\mathrm{d}x<+\infty.\notag
\end{equation*}
Moreover, if $\alpha\leq2e^{-\left(\inf\limits_{s\geq1}\chi(s)\right)^{-\frac{1}{2}}}$, then
\begin{equation*}
  \sup_{u\in F_1,\|u\|_{\omega_1}\leq1} \int_{\mathbb{R}^2}\left(e^{\alpha\left(e^{2\pi u^2}-1\right)}-1\right)\mathrm{d}x<+\infty.
\end{equation*}
If $\alpha>2 \, e^{\int^{+\infty}_{0}\log^2(1+t)e^{-2t}\mathrm{d}t}$, then
\begin{equation*}
  \sup_{u\in F_1,\|u\|_{\omega_1}\leq1} \int_{\mathbb{R}^2}\left(e^{\alpha\left(e^{2\pi u^2}-1\right)}-1\right)\mathrm{d}x=+\infty.
\end{equation*}

The similar results in $\mathbb{R}^N(N\geq2)$ was provided by Aouaoui and Jlel \cite{SRn}, considering the following weight:
\begin{equation}\label{omega2}
  \omega_2(x)=\begin{cases}
  \big(\log(\frac{e}{\mid x\mid})\big)^{\beta(N-1)},&\ \text{if $|x|<1$},\\
  \chi_1(|x|),&\ \text{if $|x|\geq1$},
  \end{cases}
\end{equation}
where $0<\beta\leq1$ and $\chi_1:[1,+\infty)\rightarrow(0,+\infty)$ is a continuous function such that $\chi_1(1)=1$ and $\inf\limits_{t\geq1}\chi_1(t)>0$. Denote by $E_\beta$ as the space of all radial functions of the completion of $C^{\infty}_0(\mathbb{R}^N)$ with respect to the norm
\begin{equation*}
  \|u\|^N_{\omega_2}=\int_{\mathbb{R}^N}|\nabla u|^N\omega_2(x)\mathrm{d}x.
\end{equation*}
More details see \cite{SRn}, here we just mention the main result of the case $0<\beta<1$ as follows.

\noindent\textbf{Theorem B}
 Let $0<\beta<1$ and $\omega_2$ be defined by (\ref{omega2}). Assume that there exists a positive $C>0$ such that
\begin{equation*}
  \int^{+\infty}_{1}\left(\int^{+\infty}_{r}\frac{\mathrm{d}s}{s(\chi_1(s))^{\frac{1}{N-1}}}\right)^{\frac{k}{1-\beta}}r^{N-1}\mathrm{d}r
  \leq C^k, \quad\text{for all}\,\,\, k\in \mathbb{N}^*.
\end{equation*}
For all $\alpha>0$ and $u\in E_\beta$, we have
\begin{equation*}
  \int_{\mathbb{R}^N}\left(e^{\alpha|u|^{\frac{N}{(1-\beta)(N-1)}}}-1\right)\mathrm{d}x<+\infty.\notag
\end{equation*}
Moreover,
\begin{equation*}
  \sup_{u\in E_\beta,\|u\|_{\omega_2}\leq1} \int_{\mathbb{R}^N}\left(e^{\alpha|u|^{\frac{N}{(1-\beta)(N-1)}}}-1\right)\mathrm{d}x<+\infty,
\end{equation*}
if and only if
\begin{equation*}
  \alpha\leq\alpha_{N,\beta}=N\left[\omega_{N-1}^{\frac{1}{N-1}}(1-\beta)\right]^{\frac{1}{1-\beta}}.
\end{equation*}

Extension of the Trudinger-Moser inequality to higher order Sobolev spaces can see \cite{ZW}. Zhu and Wang \cite{ZW} proved the following Adams' inequality when the domain is a unit ball in $\mathbb{R}^4$.

\noindent\textbf{Theorem C}
Let $B=B_1(0)\subset \mathbb{R}^4$ denote the unit ball in $\mathbb{R}^4$ and $\beta\in(0,1)$.
Then
\begin{equation}\label{thrm1.1}
  \sup\limits_{u\in W^{2,2}_{0,rad}(B,\omega),{\|\triangle u\|}_{\omega}\leq1}
  \int^{}_{B}e^{\alpha|u|^{\frac{2}{1-\beta}}}\mathrm{d}x<+\infty,
\end{equation}
if and only if
\begin{equation*}
  \alpha\leq\alpha_\beta=4\left[8\pi^2(1-\beta)\right]^{\frac{1}{1-\beta}},
\end{equation*}
where
$\omega(x)=\big(\log\frac{e}{|x|}\big)^{\beta}$, and
\begin{equation*}
    W^{2,2}_{0,rad}(B,\omega)=cl\Big\{u\in C^{\infty}_{0,rad}(B);\| \triangle u\|^2_{\omega}=\int^{}_{B}|\triangle u|^2\omega(x)\mathrm{d}x<+\infty\Big\}.
\end{equation*}

A nature question is whether the similar result as (\ref{thrm1.1}) holds on unbounded domain. In this paper, we focus on the bi-Laplacian operator in $\mathbb{R}^4$. We define the radial weight
\begin{equation}\label{omega}
  \omega_\beta(x)=\begin{cases}
  \big(\log(\frac{e}{\mid x\mid})\big)^{\beta},&\ \text{if $|x|<1$},\\
  \chi(|x|),&\ \text{if $|x|\geq1$},
  \end{cases}
\end{equation}
where $0<\beta<1$, $\chi:[1,+\infty)\rightarrow(0,+\infty)$ is a continuous function such that $\chi(1)=1$, $\inf\limits_{t\geq1}\chi(t)>0$ and satisfies
\begin{equation}\label{chi0}
  \int^{+\infty}_{1}\chi(y)^{-1}y^{-1}\mathrm{d}y<+\infty,
\end{equation}
and
\begin{equation}\label{chi0(1)}
  \int^{+\infty}_{1}\chi(y)^{-1}y\mathrm{d}y<+\infty.
\end{equation}

We define the weighted Sobolev space
\begin{equation*}
  X_\beta=\Big\{u\in L^{\frac{2}{1-\beta}}_{rad}(\mathbb{R}^4),\int^{}_{R^4}\omega_{\beta}(x)|\triangle u|^2\mathrm{d}x<+\infty\Big\},\quad\text{for all}\,\,\,0<\beta<1,
\end{equation*}
with the norm
\begin{equation*}
 \|u\|^2_{\beta}=\int^{}_{\mathbb{R}^4}\omega_{\beta}(x)|\triangle u|^2\mathrm{d}x.
\end{equation*}

For the function $\chi$, we make the following assumptions:
\noindent there exists a constant $M>0$ such that
\begin{equation}\label{chi1}
  \frac{1}{r^8}\left(\int^{r}_{1}t^3\chi(t)\mathrm{d}t\right)\left(\int^{r}_{1}\frac{t^3}{\chi(t)}\mathrm{d}t\right)\leq M,\quad\text{for all} \,\,\, r\geq1;
\end{equation}
\begin{equation}\label{chi2}
  \frac{1}{r^8}\int^{r}_{1}t^3\chi(t)\mathrm{d}t\leq M,\quad\text{for all} \,\,\,r\geq1;
\end{equation}
\begin{equation}\label{chi3}
  \frac{\max\limits_{r\leq t\leq4r}\chi(t)}{\min\limits_{r\leq t\leq4r}\chi(t)}\leq M,\quad\text{for all}\,\,\, r\geq1.
\end{equation}
The conditions (\ref{chi1}), (\ref{chi2}) and (\ref{chi3}) are used to guarantee the weight $\omega_\beta$ is in Muckenhoupt's $A_2$-class, that is,
$\omega_\beta$ satisfies
\begin{equation}\label{A2weight}
   \sup_{B\subset\mathbb{R}^4}\left(\frac{1}{|B|}\int^{}_{B}\omega_{\beta}(x)\mathrm{d}x\right)
  \left(\frac{1}{|B|}\int^{}_{B}\omega_{\beta}(x)^{-1}\mathrm{d}x\right)<+\infty,
\end{equation}
see \cite{TK,AC,ENK} for more details.

The main results are as follows.
\begin{theorem}\label{conclu1}
Let $0<\beta<1$ and $\omega_{\beta}$ be defined by (\ref{omega}). We suppose that

$(D_{1})$
\begin{equation}\label{D1}
  \int^{+\infty}_{1}r^3\left(\int^{+\infty}_{r}\frac{1}{\chi(y)y}\mathrm{d}y\right)^\frac{k}{1-\beta}\mathrm{d}r<C_1^k, \, C_1>0.\notag
\end{equation}

$(D_{2})$
\begin{equation}\label{D2}
  \int^{+\infty}_{1}r^{3-\frac{2k}{1-\beta}}\left(\int^{r}_{1}\frac{y}{\chi(y)}\mathrm{d}y\right)^\frac{k}{1-\beta}\mathrm{d}r<C_2^k, \, C_2>0.\notag
\end{equation}
For $u\in X_\beta$, we have
\begin{equation}\label{therom1.2(1)}
  \int^{}_{\mathbb{R}^4}\left(e^{|u|^{\frac{2}{1-\beta}}}-1\right)\mathrm{d}x<+\infty.
\end{equation}
Moreover, if $\alpha\leq\alpha_\beta$, then
\begin{equation}\label{therom1.2(2)}
  \sup_{u\in X_\beta,\|u\|_\beta\leq1}\int^{}_{\mathbb{R}^4}\left(e^{\alpha|u|^{\frac{2}{1-\beta}}}-1\right)\mathrm{d}x<+\infty.
\end{equation}
If $\alpha>\alpha_\beta$, then
\begin{equation*}
  \sup_{u\in X_\beta,\|u\|_\beta\leq1}\int^{}_{\mathbb{R}^4}\left(e^{\alpha|u|^{\frac{2}{1-\beta}}}-1\right)\mathrm{d}x=+\infty.
\end{equation*}
\end{theorem}

\begin{remark}
{\rm An example of function $\chi$ satisfying $(\ref{chi0})$, $(\ref{chi0(1)})$, $(\ref{chi1})$, $(\ref{chi2})$, $(\ref{chi3})$, $(D_1)$ and $(D_2)$ is given by}
\begin{equation*}
  \chi(y)=y^\alpha, \,\, \max\{2,4(1-\beta)\}<\alpha<4.
\end{equation*}
\end{remark}

Notice that there is already the result of Adams' inequality in the $B\subset\mathbb{R}^4$, that is Theorem C. The key point and difficulty of the proof in Theorem \ref{conclu1} is that we need a radial lemma to deal with $u(x)$ if $|x|$ large enough. Inspired by \cite{ZW,T}, we make a suitable change of variable. The more details see Lemma \ref{radial}.
Then, we establish an Lions \cite{LPL} type result about an improvement of (\ref{therom1.2(2)}) for weakly convergent sequences in $E_\beta$ with a constant larger than that found in (\ref{therom1.2(2)}) which is used to treat nonlinear equations with exponential growth at infinity. More precisely, we prove the theorem as follows.

\begin{theorem}\label{conclu2}
Let $0<\beta<1$. Assume that $(D_1)$ and $(D_2)$ hold. Let $\{u_n\}\subset X_\beta$ and $u\in X_\beta\setminus\{0\}$ be such that
${\|u_n\|}_\beta=1$ and $u_n\rightharpoonup u$ weakly in $X_\beta$.
Then
\begin{equation}\label{therom2.1}
  \sup_{n}\int^{}_{\mathbb{R}^4}\left(e^{p\alpha_\beta{\mid u_n\mid}^\frac{2}{1-\beta}}-1\right)\mathrm{d}x<+\infty,\quad\text{for all} \,\,\, 0<p<P_\beta(u),
\end{equation}
where
\begin{equation}\label{therom2.1(2)}
  P_\beta(u)=\begin{cases}
  \frac{1}{\left(1-{\parallel u\parallel}_\beta^2\right)^\frac{1}{1-\beta}},&\ if \, {\|u\|}_\beta<1,\\
  +\infty,&\ if \,{\|u\|}_\beta=1.
  \end{cases}
\end{equation}
\end{theorem}

Next, we will use the Adams' inequalities established in Theorem \ref{conclu1} and the improved inequality given in Theorem \ref{conclu2} to study the weighted elliptic problem involving nonlinearities which satisfy some conditions including the exponential growth. For $0<\beta<1$, we consider the problem
\begin{equation}\label{problem}
  \triangle\big(\omega_\beta(x)\triangle u\big)=f(u) \quad\quad \text{in} \,\, \mathbb{R}^4.
\end{equation}

We assume

$(F_1)$ $f:\mathbb{R}\rightarrow \mathbb{R}$ is a continuous function such that $f(s)=0$, $\text{for all} \, s\leq0$ and
\begin{equation}\label{F_1}
  |f(s)|\leq c_0\left(|s|^{p-1}+|s|^{q-1}\left(e^{\alpha\mid s\mid^{\frac{2}{1-\beta}}}-1\right)\right),\quad\text{for all} \,\,\, s\in \mathbb{R},\notag
\end{equation}
with $c_0>0$, $\min\{p,q\}>\frac{2}{1-\beta}$, $\alpha>0$.

$(F_2)$ There exist $M_0>0$ and $s_0>0$ such that
\begin{equation}\label{F_3}
  F(s)=\int^{s}_{0}f(t)dt\leq M_{0}f(s),\quad\text{for all} \,\,\, s\geq s_0.\notag
\end{equation}

$(F_3)$ There exists $\kappa>\frac{2}{\pi^2e^4}\left(\frac{\alpha_\beta}{\alpha}\right)^{1-\beta}$ such that $\liminf\limits_{s\rightarrow\infty}\frac{f(s)s}{e^{\alpha s^{\frac{2}{1-\beta}}}}\geq \kappa$.

\begin{theorem}\label{conclu3}
Let $0<\beta<1$. Assume that $(D_1)$, $(D_2)$ and $(F_1)-(F_3)$ hold. Then, problem (\ref{problem}) has at least one nontrivial and nonnegative weak solution.
\end{theorem}

\begin{remark}
{\rm An example of function $f$ satisfying the conditions $(F_1)-(F_3)$ is given by
\begin{equation*}
  f(s)=\begin{cases}
  Cs^{p-1}\left(\frac{2\alpha}{1-\beta}s^{\frac{2}{1-\beta}}+1\right)\left(e^{\alpha s^{\frac{2}{1-\beta}}}-1\right),&\quad \text{if} \,\, s\geq0,\\
  0,&\quad \text{if} \,\, s<0,
  \end{cases}
\end{equation*}
where $p>\frac{2}{1-\beta}$, $\alpha>0$ and $C>0$ is chosen large enough.

Indeed, we have
\begin{equation*}
  F(s)=\int^{s}_{0}f(t)\mathrm{d}t=\begin{cases}
  Cs^p\left(e^{\alpha s^{\frac{2}{1-\beta}}}-1\right),&\quad \text{if} \,\, s\geq0,\\
  0,&\quad \text{if} \,\, s<0.
  \end{cases}
\end{equation*}}
\end{remark}

\begin{remark}
{\rm In general, we often use the following (AR) condition to guarantee the Mountain Pass geometry and the boundedness of the (PS) sequences.

\noindent (AR) There exists $\theta>2$ such that $0<\theta F(t)\leq f(t)t$ for all $t\geq0$.

In this paper, we only use the conditions $(F_1)$ and $(F_3)$ to prove the Mountain Pass geometry of $I_\beta$ and use the condition $(F_2)$ to guarantee the boundedness of (PS) sequences. The main difficulty of the proof is the lack of compactness. In fact, we can prove that the (PS) condition is satisfied at a certain level of energy using Theorem \ref{conclu2} and $(F_1)$ mainly. Inspired by \cite{YYYY}, by considering the family of Adams type function, we will estimate the level of the energy $I_\beta$ using the condition $(F_3)$. It is worth mentioning that the condition $(F_3)$ we give in this paper is weaker than the corresponding condition $(H_5)$ in \cite{YYYY}. The more details see Section 5.}
\end{remark}

\begin{remark}
{\rm The conditions (\ref{chi0}) and (\ref{chi0(1)}) do not work in the proof of Theorem \ref{conclu1} and Theorem \ref{conclu2}, but are crucial to the proof of Theorem \ref{conclu3}. More specifically, the conditions play an important part in Proposition \ref{embedding}. Indeed, the condition (\ref{chi0}) is to guarantee the boundedness of $u(e_1)$. The condition (\ref{chi0(1)}) is to guarantee $u(x)\rightarrow0$ when $|x|\rightarrow\infty$ which is sufficient to prove the compactness of the embedding. More details see Proposition \ref{embedding}.}
\end{remark}

This paper is organized as follows: In Section 2, we prove some results which will be used in the next sections. Theorem \ref{conclu1} will be proved in Section 3. Section 4 is devoted for a concentration compactness result of Lions type which will be used to prove the application in Section 5. The constant $C$ in this paper will change according to different situations.

\section{Preliminaries}
In this section, we establish an estimate of $u(x)$ for $|x|$ large enough which plays an important role in the proof of Theorem \ref{conclu1}. First, we prove the following density property.
\begin{proposition}\label{prop1}
If $\chi$ satisfies the conditions (\ref{chi1}), (\ref{chi2}) and (\ref{chi3}), then $C^{\infty}_{0,rad}(\mathbb{R}^4)$ is dense in $X_\beta$.
\end{proposition}
\begin{proof}
Similar to \cite{SR,SRn}, we only need to prove $\omega_{\beta}\in A_2$, that is, we show that (\ref{A2weight}) holds.

\noindent Let $r>0$ and $x_0\in \mathbb{R}^4$. Denote by $B(x_0,r)$ the open ball of $\mathbb{R}^4$ of center $x_0$ and radius $r$.

\noindent\emph {Case 1.} $B(x_0,r)\bigcap B(0,r)\neq\phi$. In this case $B(x_0,r)\subset B(0,3r)$ which implies that
\begin{align}\label{case1.1}
    & \frac{1}{|B(x_0,r)|}\left(\int^{}_{B(x_0,r)}\omega_{\beta}(x)\mathrm{d}x\right)
      \left(\frac{1}{|B(x_0,r)|}\int^{}_{B(x_0,r)}\omega_{\beta}(x)^{-1}\mathrm{d}x\right)\notag \\
    & \leq \frac{c}{r^8}\left(\int^{3r}_{0}\omega_{\beta}(t)t^3\mathrm{d}t\right)\left(\int^{3r}_{0}\frac{t^3}{\omega_{\beta}(t)}\mathrm{d}t\right).
\end{align}

If $3r<1$, after a simple calculation, we have
\begin{equation}\label{case1.2}
  \limsup_{r\rightarrow0^+}\frac{c}{r^8}\left(\int^{3r}_{0}t^3(1-\log t)^{\beta}\mathrm{d}t\right)
  \left(\int^{3r}_{0}\frac{t^3}{(1-\log t)^{\beta}}\mathrm{d}t\right)<+\infty.
\end{equation}

If $3r\geq1$, then
\begin{align}\label{case1.3}
    & \frac{c}{r^8}\left(\int^{3r}_{0}\omega_{\beta}(t)t^3\mathrm{d}t\right)\left(\int^{3r}_{0}\frac{t^3}{\omega_{\beta}(t)}\mathrm{d}t\right)\notag \\
    & =\frac{c}{r^8}\left(\int^{1}_{0}t^3(1-\log t)^{\beta}\mathrm{d}t+\int^{3r}_{1}t^3\chi(t)\mathrm{d}t\right)\notag \\
    & \times\left(\int^{1}_{0}\frac{t^3}{(1-\log t)^{\beta}}\mathrm{d}t+\int^{3r}_{1}\frac{t^3}{\chi(t)}\mathrm{d}t\right).
\end{align}
Since $\inf\limits_{t\geq1}\chi(t)>0$, then
\begin{equation}\label{case1.4}
  \limsup_{r\rightarrow+\infty}\frac{1}{r^8}\left(\int^{3r}_{1}\frac{t^3}{\chi(t)}\mathrm{d}t\right)=0<+\infty.
\end{equation}
Combining (\ref{chi1}), (\ref{chi2}), (\ref{case1.2}), (\ref{case1.3}) and (\ref{case1.4}), we deduce from (\ref{case1.1}) that
\begin{equation}\label{case1.5}
  \frac{1}{|B(x_0,r)|}\left(\int^{}_{B(x_0,r)}\omega_{\beta}(x)\mathrm{d}x\right)
  \left(\frac{1}{|B(x_0,r)|}\int^{}_{B(x_0,r)}\omega_{\beta}(x)^{-1}\mathrm{d}x\right)<\infty.
\end{equation}
\emph {Case 2.} $B(x_0,r)\bigcap B(0,r)=\phi$. In this case, we have
\begin{equation*}
  \frac{|x_0|}{2}\leq |x|\leq 2|x_0|,\quad\text{for all} \,\,\, x\in B(x_0,r).
\end{equation*}
Hence,
\begin{align}\label{case2.1}
    & \left(\frac{1}{|B(x_0,r)|}\int^{}_{B(x_0,r)}\omega_{\beta}(x)\mathrm{d}x \right)
  \left(\frac{1}{|B(x_0,r)|}\int^{}_{B(x_0,r)}{\omega_{\beta}(x)}^{-1}\mathrm{d}x \right)\notag \\
    & \leq\frac{\sup\limits_{\frac{\mid x_0\mid}{2}\leq\mid x\mid\leq2\mid x_0\mid}\omega_{\beta}(x)}
  {\inf\limits_{\frac{|x_0|}{2}\leq |x| \leq 2|x_0|}\omega_{\beta}(x)}\notag \\
    & \leq\sup\limits_{\alpha>0} \left(\frac{\sup\limits_{\alpha\leq t\leq 4\alpha}\omega_{\beta}(t)}
  {\inf\limits_{\alpha\leq t\leq 4\alpha}\omega_{\beta}(t)}\right).
\end{align}

If $4\alpha<1$, then
\begin{equation*}
  \frac{\sup\limits_{\alpha\leq t\leq 4\alpha}\omega_{\beta}(t)}{\inf\limits_{\alpha\leq t\leq 4\alpha}\omega_{\beta}(t)}
  =\frac{(1-\log\alpha)^{\beta}}{(1-\log(4\alpha))^{\beta}},
\end{equation*}
and consequently,
\begin{equation}\label{case2.2}
  \sup\limits_{0<\alpha<\frac{1}{4}} \left(\frac{\sup\limits_{\alpha\leq t\leq 4\alpha}\omega_{\beta}(t)}{\inf\limits_{\alpha\leq t\leq 4\alpha}\omega_{\beta}(t)}\right)
  =\sup\limits_{0<\alpha<\frac{1}{4}} \left(\frac{1-\log\alpha}{1-\log(4\alpha)}\right)^{\beta} <+\infty.
\end{equation}

If $\frac{1}{4}\leq\alpha<1$, by the continuity of the function $\omega_{\beta}$ and the boundedness theorem, we get
\begin{equation*}
  \frac{\sup\limits_{\alpha\leq t\leq 4\alpha}\omega_{\beta}(t)}{\inf\limits_{\alpha\leq t\leq 4\alpha}\omega_{\beta}(t)}
  \leq \frac{\sup\limits_{\frac{1}{4}\leq t\leq4}\omega_{\beta}(t)}{\inf\limits_{\frac{1}{4}\leq t\leq4}\omega_{\beta}(t)} <+\infty,
\end{equation*}
and consequently,
\begin{equation}\label{case2.3}
  \sup\limits_{\frac{1}{4}\leq\alpha<1}\left(\frac{\sup\limits_{\alpha\leq t\leq 4\alpha}\omega_{\beta}(t)}
  {\inf\limits_{\alpha\leq t\leq 4\alpha}\omega_{\beta}(t)}\right)
  <+\infty.
\end{equation}

If $\alpha\geq1$, then by (\ref{chi3}) we have
\begin{equation*}
  \frac{\sup\limits_{\alpha\leq t\leq 4\alpha}\omega_{\beta}(t)}{\inf\limits_{\alpha\leq t\leq 4\alpha}\omega_{\beta}(t)}
  \leq \frac{\sup\limits_{\alpha\leq t\leq 4\alpha}\chi(t)}{\inf\limits_{\alpha\leq t\leq 4\alpha}\chi(t)} <+\infty,
\end{equation*}
and consequently,
\begin{equation}\label{case2.4}
  \sup\limits_{\alpha\geq1}\left(\frac{\sup\limits_{\alpha\leq t\leq 4\alpha}\omega_{\beta}(t)}{\inf\limits_{\alpha\leq t\leq 4\alpha}\omega_{\beta}(t)}\right) <+\infty.
\end{equation}
Combining (\ref{case2.2}), (\ref{case2.3}) and (\ref{case2.4}), we deduce from (\ref{case2.1}) that
\begin{equation}\label{case2.5}
  \frac{1}{|B(x_0,r)|}\left(\int^{}_{B(x_0,r)}\omega_{\beta}(x)\mathrm{d}x\right)
  \left(\frac{1}{|B(x_0,r)|}\int^{}_{B(x_0,r)}\omega_{\beta}(x)^{-1}\mathrm{d}x\right)<\infty.
\end{equation}

From (\ref{case1.5}) and (\ref{case2.5}), we have $\omega_{\beta}\in A_2$. Then, $C^{\infty}_{0,rad}(R^4)$ is dense in $X_\beta$.
\end{proof}

Now, using the density of $C^{\infty}_{0,rad}(\mathbb{R}^4)$ in $X_\beta$, we can prove the following radial lemma which is based on a suitable change of variable, see \cite{ZW,T}.
\begin{lemma}\label{radial}
Fix $0<\beta<1$ and let $u\in X_{\beta}$. Then one has

\noindent

\begin{equation}\label{radials}
 \begin{split}
  & |u(x)|\leq \frac{\sqrt{2}c_1}{4\pi} \left(\int^{}_{|x|\geq1}|\triangle u|^2\chi(|x|)\mathrm{d}x \right)^{\frac{1}{2}}
  \left[\left( \int^{+\infty}_{|x|} \chi(y)^{-1}y^{-1} \mathrm{d}y\right)^{\frac{1}{2}}
  +\left(\int^{|x|}_{1} \chi(y)^{-1}y|x|^{-2} \mathrm{d}y \right)^{\frac{1}{2}}\right]+\frac{c_2}{|x|^2}, \\
  & for  \,\, |x|\geq1.
 \end{split}
\end{equation}
\end{lemma}
\begin{proof}
Since $C^{\infty}_{0,rad}(\mathbb{R}^4)$ is dense in $X_\beta$, it is enough to prove (\ref{radials}) when $u\in C^{\infty}_{0,rad}(\mathbb{R}^4)$.
Setting
\begin{equation*}
  \omega(s)=4\pi u(x),\; r=|x|=s^{-\frac{1}{2}},\quad\text{for} \,\,\, 0<s\leq1.
\end{equation*}
Then we have
\begin{align*}
  &\lim_{s\rightarrow0^{+}}\omega(0)=0,\lim_{s\rightarrow1^{-}} \omega'(s)=\omega^{'}(1)=-2\pi u'(1), \; c_2\triangleq\left|-2\pi u'(1)\right|\\
  &\omega''(s)=\pi \cdot s^{-3}\triangle u\mid_{|x|=s^{-\frac{1}{2}}},
\end{align*}
\noindent and
\begin{equation*}
  \begin{split}
  \int^{}_{|x|\geq1}|\triangle u|^2 \chi(|x|)\mathrm{d}x & =2\pi^2 \int^{+\infty}_{1}|\triangle u|^2 \chi(r)r^3\mathrm{d}r \\
  & =\pi^2\int^{1}_{0}\left|\triangle u\left(s^{-\frac{1}{2}}\right)\right|^2 \chi\left(s^{-\frac{1}{2}}\right)s^{-3}\mathrm{d}s \\
  & =\int^{1}_{0}\left|\omega''(s)\right|^2 \chi\left(s^{-\frac{1}{2}}\right)s^3\mathrm{d}s.
  \end{split}
\end{equation*}
Let us now evaluate $\omega(s)$, we have
\begin{equation*}
   \begin{split}
     |\omega(s)| & =\left|\int^{s}_{0}\omega'(z)\mathrm{d}z\right| \\
     & =\left|\int^{s}_{0} \left(\int^{1}_{z}\omega''(t)\mathrm{d}t-\omega^{'}(1)\right) \mathrm{d}z\right| \\
     & \leq \int^{s}_{0}\left| \int^{1}_{z}\omega''(t)\mathrm{d}t \right|\mathrm{d}z+\left| \int^{s}_{0}c_2 \mathrm{d}z \right| \\
     & \leq \left( \int^{s}_{0}\int^{1}_{z} |\omega''(t)|^2\chi\left(t^{-\frac{1}{2}}\right)t^2 \mathrm{d}t\mathrm{d}z \right)^{\frac{1}{2}}
       \cdot \left( \int^{s}_{0}\int^{1}_{z} \chi\left(t^{-\frac{1}{2}}\right)^{-1}t^{-2} \mathrm{d}t\mathrm{d}z \right)^{\frac{1}{2}}+|c_2s| \\
     & \leq \left( \int^{s}_{0} |\omega''(t)|^2\chi\left(t^{-\frac{1}{2}}\right)t^2 \mathrm{d}t \int^{t}_{0}\mathrm{d}z
       +\int^{1}_{s} |\omega''(t)|^2\chi\left(t^{-\frac{1}{2}}\right)t^2 \mathrm{d}t \int^{s}_{0}\mathrm{d}z \right)^{\frac{1}{2}} \\
       & \cdot \left( \int^{s}_{0}\int^{1}_{z} \chi\left(t^{-\frac{1}{2}}\right)^{-1}t^{-2} \mathrm{d}t\mathrm{d}z \right)^{\frac{1}{2}}+|c_2s| \\
     & \leq \left( \int^{s}_{0} |\omega''(t)|^2\chi\left(t^{-\frac{1}{2}}\right)t^3 \mathrm{d}t
       +\int^{1}_{s} |\omega''(t)|^2\chi\left(t^{-\frac{1}{2}}\right)t^3 \mathrm{d}t \right)^{\frac{1}{2}} \\
       & \cdot \left( \int^{s}_{0}\int^{1}_{z} \chi\left(t^{-\frac{1}{2}}\right)^{-1}t^{-2} \mathrm{d}t\mathrm{d}z \right)^{\frac{1}{2}}+|c_2s| \\
     & =\left( \int^{1}_{0} |\omega''(t)|^2\chi\left(t^{-\frac{1}{2}}\right)t^3 \mathrm{d}t \right)^{\frac{1}{2}}
       \left(\int^{s}_{0}\int^{1}_{z} \chi\left(t^{-\frac{1}{2}}\right)^{-1}t^{-2} \mathrm{d}t\mathrm{d}z \right)^{\frac{1}{2}}+|c_2s| \\
     & =\left( \int^{}_{|x|\geq1} |\triangle u|^2\chi\left(|x|\right) \mathrm{d}x \right)^{\frac{1}{2}}
       \left(\int^{s}_{0}\int^{1}_{z} \chi\left(t^{-\frac{1}{2}}\right)^{-1}t^{-2} \mathrm{d}t\mathrm{d}z \right)^{\frac{1}{2}}+|c_2s|,
   \end{split}
\end{equation*}
and
\begin{equation*}
  \begin{split}
     \int^{s}_{0}\int^{1}_{z} \chi\left(t^{-\frac{1}{2}}\right)^{-1}t^{-2} \mathrm{d}t\mathrm{d}z
     & =\int^{s}_{0} \chi\left(t^{-\frac{1}{2}}\right)^{-1}t^{-2}\mathrm{d}t \int^{t}_{0}\mathrm{d}z
      +\int^{1}_{s} \chi\left(t^{-\frac{1}{2}}\right)^{-1}t^{-2}\mathrm{d}t \int^{s}_{0}\mathrm{d}z \\
     & =\int^{s}_{0} \chi\left(t^{-\frac{1}{2}}\right)^{-1}t^{-1}\mathrm{d}t
      +\int^{1}_{s} \chi\left(t^{-\frac{1}{2}}\right)^{-1}t^{-2}s\mathrm{d}t.
   \end{split}
\end{equation*}
Thus we have
\begin{equation*}
  |\omega(s)|\leq \left( \int^{}_{|x|\geq1} |\triangle u|^2\chi(|x|) \mathrm{d}x \right)^{\frac{1}{2}}
  \left(\int^{s}_{0} \chi\left(t^{-\frac{1}{2}}\right)^{-1}t^{-1}\mathrm{d}t +\int^{1}_{s} \chi\left(t^{-\frac{1}{2}}\right)^{-1}t^{-2}s\mathrm{d}t \right)^{\frac{1}{2}}+|c_2s|.
\end{equation*}
Consequently,
\begin{equation*}
  \begin{split}
   |u(x)| & \leq \frac{c_1}{4\pi} \left(\int^{}_{|x|\geq1}|\triangle u|^2\chi(|x|)\mathrm{d}x \right)^{\frac{1}{2}}
            \left[\left( \int^{\frac{1}{|x|^2}}_{0} \chi\left(t^{-\frac{1}{2}}\right)^{-1}t^{-1} \mathrm{d}t\right)^{\frac{1}{2}}
            +\left( \int^{1}_{\frac{1}{|x|^2}} \chi\left(t^{-\frac{1}{2}}\right)^{-1}t^{-2}|x|^{-2} \mathrm{d}t \right)^{\frac{1}{2}}\right] \\
          & +\frac{c_2}{|x|^2} \\
          & = \frac{\sqrt{2}c_1}{4\pi} \left(\int^{}_{|x|\geq1}|\triangle u|^2\chi(|x|)\mathrm{d}x \right)^{\frac{1}{2}}
              \left[\left( \int^{+\infty}_{|x|} \chi(y)^{-1}y^{-1} \mathrm{d}y\right)^{\frac{1}{2}}
              +\left( \int^{|x|}_{1} \chi(y)^{-1}y|x|^{-2} \mathrm{d}y \right)^{\frac{1}{2}}\right]
              +\frac{c_2}{|x|^2}.
  \end{split}
\end{equation*}
\end{proof}

\section{Proof of Theorem \ref{conclu1}}
In this section, we prove the capital Adams' inequality in detail.
We start by proving the first conclusion of Theorem \ref{conclu1}. For $u\in X_\beta$, we have
\begin{equation}\label{2.1}
  \int^{}_{\mathbb{R}^4}\left(e^{|u|^{\frac{2}{1-\beta}}}-1\right)\mathrm{d}x
  =\int^{}_{|x|\geq1}\left(e^{|u|^{\frac{2}{1-\beta}}}-1\right)\mathrm{d}x
  +\int^{}_{|x|<1}\left(e^{|u|^{\frac{2}{1-\beta}}}-1\right)\mathrm{d}x.
\end{equation}

On the one hand,
\begin{equation}\label{2.2}
  \int^{}_{|x|\geq1}\left(e^{|u|^{\frac{2}{1-\beta}}}-1\right)\mathrm{d}x
  =\sum^{+\infty}_{k=1} \frac{1}{k!}\int^{}_{|x|\geq1}|u|^{\frac{2k}{1-\beta}}\mathrm{d}x.
\end{equation}
From (\ref{radials}) and $(D_1)$, $(D_2)$, we get
\begin{equation}\label{2.3}
  \begin{split}
     \int^{}_{|x|\geq1}|u|^{\frac{2k}{1-\beta}}\mathrm{d}x
     & \leq \int^{}_{|x|\geq1} \left\{\frac{\sqrt{2}c_1}{4\pi}\|u\|_\beta
      \left[\left( \int^{+\infty}_{|x|} \chi(y)^{-1}y^{-1} \mathrm{d}y\right)^{\frac{1}{2}}
      +\left( \int^{|x|}_{1} \chi(y)^{-1}y|x|^{-2} \mathrm{d}y \right)^{\frac{1}{2}}\right]
         +\frac{c_2}{|x|^2} \right\}^{\frac{2k}{1-\beta}}\mathrm{d}x \\
     & \leq C\|u\|^{\frac{2k}{1-\beta}}_\beta  \left[\int^{}_{|x|\geq1} \left(\int^{+\infty}_{|x|} \chi(y)^{-1}y^{-1}\mathrm{d}y\right)^{\frac{k}{1-\beta}}\mathrm{d}x
         +\int^{}_{|x|\geq1} \left(\int^{|x|}_{1} \chi(y)^{-1}y|x|^{-2}\mathrm{d}y\right)^{\frac{k}{1-\beta}}\mathrm{d}x \right] \\
     & +C\int^{}_{|x|\geq1}\frac{1}{|x|^{\frac{4k}{1-\beta}}}\mathrm{d}x \\
     & = C\|u\|^{\frac{2k}{1-\beta}}_\beta 2\pi^2 \left[\int^{+\infty}_{1}r^3\left(\int^{+\infty}_{r}\chi(y)^{-1}y^{-1}\right)^{\frac{k}{1-\beta}}\mathrm{d}r
       +\int^{+\infty}_{1}r^{3-\frac{2k}{1-\beta}}\left(\int^{r}_{1}\frac{y}{\chi(y)}\mathrm{d}y\right)^{\frac{k}{1-\beta}}\mathrm{d}r \right] \\
     & +2\pi^2\int^{+\infty}_{1}\frac{1}{r^{\frac{4k}{1-\beta}-3}}\mathrm{d}r <C^k        ,\quad\text{for all} \,\,\, k\geq1.
  \end{split}
\end{equation}
Combining (\ref{2.2}) and (\ref{2.3}), we have
\begin{equation}\label{2.4}
  \int^{}_{|x|\geq1}\left(e^{|u|^{\frac{2}{1-\beta}}}-1\right)\mathrm{d}x \leq\sum^{+\infty}_{k=1}\frac{C^k}{k!}=e^{C}-1<+\infty.
\end{equation}
Now in order to estimate the second integral in (\ref{2.1}), set
\begin{equation}\label{vx}
  v(x)=\begin{cases}
  u(x)-u(e_1),&0\leq|x|<1,\\
  0,&|x|\geq1,
  \end{cases}
\end{equation}
where $e_1=(1,0,0,0)\in \mathbb{R}^4$. Clearly $v\in W^{2,2}_{0,rad}(B,\omega)$.

\noindent After an elementary calculation, we have
\begin{equation*}
  |u|^{\frac{2}{1-\beta}}=|v+u(e_1)|^{\frac{2}{1-\beta}} \leq(1+\epsilon)|v|^{\frac{2}{1-\beta}}+\left(1-\frac{1}{(1+\epsilon)^{\frac{1-\beta}{1+\beta}}}\right)^{\frac{1+\beta}{\beta-1}} |u(e_1)|^{\frac{2}{1-\beta}},
\end{equation*}
for fix $\epsilon>0$.

\noindent Then, using \cite[Lemma 2(ii)]{DBJR}, we have
\begin{equation}\label{2.5}
  \begin{split}
    \int^{}_{|x|<1}e^{|u|^{\frac{2}{1-\beta}}}\mathrm{d}x
    & \leq\int^{}_{|x|<1}e^{(1+\epsilon)|v|^{\frac{2}{1-\beta}}}
      \cdot e^{\left(1-\frac{1}{(1+\epsilon)^{\frac{1-\beta}{1+\beta}}}\right)^{\frac{1+\beta}{\beta-1}} |u(e_1)|^{\frac{2}{1-\beta}}}\mathrm{d}x \\
    & \leq e^{\left(1-\frac{1}{(1+\epsilon)^{\frac{1-\beta}{1+\beta}}}\right)^{\frac{1+\beta}{\beta-1}} |u(e_1)|^{\frac{2}{1-\beta}}}  \int^{}_{|x|<1}e^{(1+\epsilon)|v|^{\frac{2}{1-\beta}}}\mathrm{d}x<+\infty.
  \end{split}
\end{equation}
Combining (\ref{2.1}), (\ref{2.4}) and (\ref{2.5}), we conclude that
\begin{equation*}
  \int^{}_{\mathbb{R}^4}\left(e^{|u|^{\frac{2}{1-\beta}}}-1\right)\mathrm{d}x<+\infty,\quad\text{for all} \,\,\, u\in X_\beta.
\end{equation*}
This ends the proof of (\ref{therom1.2(1)}).

On the other hand, for any $u\in X_\beta$, $\|u\|_\beta \leq1$, notice that
\begin{equation}\label{2.6}
  \|\triangle v\|^2_\omega=\int^{}_{B}|\triangle v|^2\omega(x)\mathrm{d}x=\int^{}_{B}|\triangle u|^2\omega_\beta(x)\mathrm{d}x\leq\|u\|^2_\beta \leq1.
\end{equation}
Let $\alpha<\alpha_\beta$. Plainly, there exists $\epsilon>0$ such that $\alpha(1+\epsilon)<\alpha_\beta$.

\noindent Combining (\ref{2.6}) and (\ref{thrm1.1}), we have
\begin{equation}\label{2.7}
 \begin{split}
     \sup_{u\in X_\beta,\|u\|_\beta\leq1}\int^{}_{|x|<1}e^{\alpha|u|^{\frac{2}{1-\beta}}}\mathrm{d}x
     & \leq C\sup_{u\in X_\beta,\|u\|_\beta\leq1} \int^{}_{|x|<1}e^{\alpha(1+\epsilon)|v|^{\frac{2}{1-\beta}}}\mathrm{d}x \\
     & \leq C\sup_{v\in W^{2,2}_{0,rad}(B,\omega),\|v\|_\beta\leq1} \int^{}_{|x|<1}e^{\alpha_\beta|v|^{\frac{2}{1-\beta}}}\mathrm{d}x<+\infty.
 \end{split}
\end{equation}
Furthermore
\begin{equation}\label{2.8}
  \int^{}_{|x|\geq1}\left(e^{\alpha|u|^{\frac{2}{1-\beta}}}-1\right)\mathrm{d}x
  =\sum^{+\infty}_{k=1} \frac{\alpha^k}{k!}\int^{}_{|x|\geq1}|u|^{\frac{2k}{1-\beta}}\mathrm{d}x.
\end{equation}
Combining (\ref{2.3}) and (\ref{2.8}), we infer
\begin{equation}\label{2.9}
  \sup_{u\in X_\beta,\|u\|_\beta\leq1}\int^{}_{|x|\geq1}\left(e^{\alpha|u|^{\frac{2}{1-\beta}}}-1\right)\mathrm{d}x<+\infty.
\end{equation}
It follows from (\ref{2.7}) and (\ref{2.9}) that
\begin{equation*}
  \sup_{u\in X_\beta,\|u\|_\beta\leq1}\int^{}_{\mathbb{R}^4}\left(e^{\alpha|u|^{\frac{2}{1-\beta}}}-1\right)\mathrm{d}x<+\infty,\quad\text{for all} \,\,\, \alpha<\alpha_\beta.
\end{equation*}

Now we prove that
\begin{equation*}
  \sup_{u\in X_\beta,\|u\|_\beta\leq1}\int^{}_{\mathbb{R}^4}\left(e^{\alpha_\beta|u|^{\frac{2}{1-\beta}}}-1\right)\mathrm{d}x<+\infty.
\end{equation*}
For any $u\in X_\beta$, $\|u\|_\beta=\int_{\mathbb{R}^4}|\triangle u|^2\omega_\beta(x)\mathrm{d}x\leq1$, we discuss it in three cases.

\noindent\emph {Case 1.} Assume that
\begin{equation*}
  \int_{|x|<1}|\triangle u|^2\omega_\beta(x)\mathrm{d}x=1.
\end{equation*}
Then
\begin{equation}\label{ab1}
  \int_{|x|\geq1}|\triangle u|^2\omega_\beta(x)\mathrm{d}x=0.
\end{equation}
From (\ref{radials}), we have
\begin{equation}\label{ab2}
  |u(x)|\leq c_2,\quad\text{for all} \,\,\, |x|\geq1.
\end{equation}
\noindent Having in mind that $u\in C^\infty_{0,rad}(\mathbb{R}^4)$.
Then combining (\ref{ab1}) and (\ref{ab2}), we deduce that for any $|x|\geq1$, $u(x)=0$. In particular, we have $u(1)=0$. So, $u\mid_B\in W^{2,2}_{0,rad}(B,\omega)$. Therefore, by Theorem C, we get
\begin{equation*}
  \int_{|x|<1}\left(e^{\alpha_\beta|u|^{\frac{2}{1-\beta}}}-1\right)\mathrm{d}x
  \leq\int_{|x|<1}e^{\alpha_\beta|u|^{\frac{2}{1-\beta}}}\mathrm{d}x
  \leq \sup\limits_{v\in W^{2,2}_{0,rad}(B,\omega),{\|\triangle v\|}_{\omega}\leq1}\int^{}_{B}e^{\alpha_\beta|v|^{\frac{2}{1-\beta}}}\mathrm{d}x<+\infty.
\end{equation*}
\emph {Case 2.} Assume that
\begin{equation*}
  \int_{|x|<1}|\triangle u|^2\omega_\beta(x)\mathrm{d}x=0.
\end{equation*}
We set $v(x)$ as (\ref{vx}), similar to (\ref{2.5}), we have
\begin{equation}\label{ab3}
  \int^{}_{|x|<1}\left(e^{\alpha_\beta|u|^{\frac{2}{1-\beta}}}-1\right)\mathrm{d}x
  \leq\int^{}_{|x|<1}e^{\alpha_\beta|u|^{\frac{2}{1-\beta}}}\mathrm{d}x\leq C\int_{|x|<1}e^{\alpha_\beta(1+\epsilon)|v|^{\frac{2}{1-\beta}}}\mathrm{d}x.
\end{equation}
Notice that $(1+\epsilon)^{\frac{1-\beta}{2}}v\in W^{2,2}_{0,rad}(B,\omega)$ and
\begin{equation*}
  \left\|\triangle\left((1+\epsilon)^{\frac{1-\beta}{2}}v\right)\right\|^2_{\omega}
  =(1+\epsilon)^{1-\beta}\int_{|x|<1}|\triangle v|^2\omega(x)\mathrm{d}x
  =(1+\epsilon)^{1-\beta}\int_{|x|<1}|\triangle u|^2\omega_\beta(x)\mathrm{d}x=0<1.
\end{equation*}
Then, by Theorem C, we get
\begin{equation*}
  \int^{}_{|x|<1}\left(e^{\alpha_\beta|u|^{\frac{2}{1-\beta}}}-1\right)\mathrm{d}x
  \leq C\int_{|x|<1}e^{\alpha_\beta\left|(1+\epsilon)^{\frac{1-\beta}{2}}v\right|^{\frac{2}{1-\beta}}}\mathrm{d}x
   \leq \sup\limits_{v\in W^{2,2}_{0,rad}(B,\omega),{\|\triangle v\|}_{\omega}\leq1}\int^{}_{B}e^{\alpha_\beta|v|^{\frac{2}{1-\beta}}}\mathrm{d}x<+\infty.
\end{equation*}
\emph {Case 3.} Assume that
\begin{equation*}
  0<\int_{|x|<1}|\triangle u|^2\omega_\beta(x)\mathrm{d}x<1.
\end{equation*}
Similar to the previous case, we also have (\ref{ab3}).

\noindent Let $1+\epsilon=\left(\int_{|x|<1}|\triangle u|^2\omega_\beta(x)\mathrm{d}x\right)^{-\frac{1}{1-\beta}}$, we have
\begin{equation*}
  \left\|\triangle\left((1+\epsilon)^{\frac{1-\beta}{2}}v\right)\right\|^2_{\omega}
  =(1+\epsilon)^{1-\beta}\int_{|x|<1}|\triangle u|^2\omega_\beta(x)\mathrm{d}x=1.
\end{equation*}
Then, by Theorem C, we get
\begin{equation*}
  \int^{}_{|x|<1}\left(e^{\alpha_\beta|u|^{\frac{2}{1-\beta}}}-1\right)\mathrm{d}x
  \leq C\int_{|x|<1}e^{\alpha_\beta\left|(1+\epsilon)^{\frac{1-\beta}{2}}v\right|^{\frac{2}{1-\beta}}}\mathrm{d}x
   \leq \sup\limits_{v\in W^{2,2}_{0,rad}(B,\omega),{\|\triangle v\|}_{\omega}\leq1}\int^{}_{B}e^{\alpha_\beta|v|^{\frac{2}{1-\beta}}}\mathrm{d}x<+\infty.
\end{equation*}
This ends the proof of (\ref{therom1.2(2)}).

Next, we show that if $\alpha>\alpha_\beta$, then the supremum is infinite.

\noindent We consider the family of Adams type function, see \cite{SFSB,ZW}
\begin{equation*}
  u_\epsilon=\begin{cases}
   \left(\frac{\log\frac{e^4}{\epsilon}}{\alpha_\beta}\right)^{\frac{1-\beta}{2}}
   -\frac{|x|^{2(1-\beta)}}{2\left(\frac{\alpha_\beta}{4}\epsilon\right)^{\frac{1-\beta}{2}}\left(\frac{1}{4}\log\frac{e^4}{\epsilon}\right)
   ^{\frac{1+\beta}{2}}}
   +\frac{1}{2\left(\frac{\alpha_\beta}{4}\right)^{\frac{1-\beta}{2}}\left(\frac{1}{4}\log\frac{e^4}{\epsilon}\right)^{\frac{1+\beta}{2}}},&|x|\leq \sqrt[4]{\epsilon},\\
   \frac{\left(\log\frac{e}{|x|}\right)^{1-\beta}} {\left(\frac{\alpha_\beta}{16}\log\frac{e^4}{\epsilon}\right)^{\frac{1-\beta}{2}}}
  ,&\sqrt[4]{\epsilon}<|x|\leq\frac{1}{2},\\
   \eta_\epsilon,&|x|>\frac{1}{2},
  \end{cases}
\end{equation*}
where $\eta\in C^\infty_0(B_1)$ is such that
\begin{equation*}
  \eta_\epsilon|_{x=\frac{1}{2}}=\frac{1}{\left(\frac{\alpha_\beta}{16}\log\frac{e^4}{\epsilon}\right)^{\frac{1-\beta}{2}}}(\log2e)^{1-\beta},
  \frac{\partial \eta_\epsilon}{\partial x}|_{x=\frac{1}{2}}
  =-\frac{2(1-\beta)}{\left(\frac{\alpha_\beta}{16}\log\frac{e^4}{\epsilon}\right)^{\frac{1-\beta}{2}}}(\log2e)^{-\beta},
\end{equation*}
and $\eta_\epsilon,\nabla\eta_\epsilon,\triangle\eta_\epsilon$ are all $O\left(\frac{1}{\left(\log\frac{e^4}{\epsilon}\right)^\frac{1-\beta}{2}}\right)$ as $\epsilon\rightarrow0^+$.

\noindent By direct computation, we get
\begin{equation*}
  \triangle u_\epsilon=\begin{cases}
  -\frac{(1-\beta)(4-2\beta)|x|^{-2\beta}}
  {\left(\frac{\alpha_\beta}{4}\epsilon\right)^{\frac{1-\beta}{2}}\left(\frac{1}{4}\log\frac{1}{\epsilon}\right)^{\frac{1+\beta}{2}}},&|x|\leq \sqrt[4]{\epsilon},\\
   \frac{1-\beta}{\left(\frac{\alpha_\beta}{16}\log\frac{1}{\epsilon}\right)^{\frac{1-\beta}{2}}}
   \left(\log\frac{1}{|x|}\right)^{-\beta}\frac{1}{|x|^2}\left(\frac{-\beta}{\log\frac{1}{|x|}}-2\right),&\sqrt[4]{\epsilon}<|x|\leq\frac{1}{2},\\
  \triangle\eta_\epsilon,&|x|>\frac{1}{2}.
  \end{cases}
\end{equation*}
Now, we evaluate the weighted Sobolev norm of $u_\epsilon$ as follows
\begin{equation*}
  \begin{split}
     \|u_\epsilon\|^2_\beta
     & =\underbrace{2\pi^2\int^{\sqrt[4]{\epsilon}}_{0} |\triangle u_\epsilon|^2\left(\log\frac{e}{|x|}\right)^\beta|x|^3\mathrm{d}|x|}_{A_1}
       +\underbrace{2\pi^2\int^{\frac{1}{2}}_{\sqrt[4]{\epsilon}} |\triangle u_\epsilon|^2\left(\log\frac{e}{|x|}\right)^\beta|x|^3\mathrm{d}|x|}_{A_2} \\
     & +\underbrace{2\pi^2\int^{1}_{\frac{1}{2}}|\triangle \eta_\epsilon|^2\left(\log\frac{e}{|x|}\right)^\beta|x|^3\mathrm{d}|x|
       +2\pi^2\int^{+\infty}_{1}|\triangle \eta_\epsilon|^2\chi(|x|)|x|^3\mathrm{d}|x|}_{A_3}.
  \end{split}
\end{equation*}
Similar to \cite{SFSB,ZW}, an easy computation gives
\begin{equation*}
 A_1=O\left(\frac{1}{\log\frac{e}{\sqrt[4]{\epsilon}}}\right),A_2=1+O\left(\frac{1}{\left(\log\frac{e}{\sqrt[4]{\epsilon}}\right)^{1-\beta}}\right),
 A_3=O\left(\frac{1}{\left(\log\frac{e}{\sqrt[4]{\epsilon}}\right)^{1-\beta}}\right).
\end{equation*}
Then
\begin{equation*}
  \|u_\epsilon\|^2_\beta\rightarrow1,\quad \text{as} \,\, \epsilon\rightarrow0^+.
\end{equation*}
Now we normalize $u_\epsilon$, set
\begin{equation*}
  \widetilde{u}_\epsilon=\frac{u_\epsilon}{\|u_\epsilon\|_\beta}\in X_\beta.
\end{equation*}
Moreover, we have
\begin{equation*}
  \widetilde{u}_\epsilon\geq\frac{\left(\frac{\log\frac{e^4}{\epsilon}}{\alpha_\beta}\right)^{\frac{1-\beta}{2}}}{\|u_\epsilon\|_\beta},\quad\text{for all} \,\,\, x\in\mathbb{R}^4, \,\, |x|\leq\sqrt[4]{\epsilon}.
\end{equation*}
So when $\alpha>\alpha_\beta$, for any $u\in X_\beta$, $\|u\|_\beta\leq1$, we have
\begin{equation*}
  \begin{split}
    \int^{}_{\mathbb{R}^4}\left(e^{\alpha|u|^{\frac{2}{1-\beta}}}-1\right)\mathrm{d}x
    & \geq\lim_{\epsilon\rightarrow0^+}\int^{}_{|x|\leq\sqrt[4]{\epsilon}}\left(e^{\alpha|\widetilde{u}_\epsilon|^{\frac{2}{1-\beta}}}-1\right)\mathrm{d}x \\
    & \geq\lim_{\epsilon\rightarrow0^+}\frac{\pi^2}{2}\epsilon\left(\left(\frac{e^4}{\epsilon}\right)^\frac{\alpha}{\alpha_\beta}-1\right)=+\infty.
  \end{split}
\end{equation*}
It follows that
\begin{equation*}
  \sup_{u\in X_\beta,\|u\|_\beta\leq1}\int^{}_{\mathbb{R}^4}\left(e^{\alpha|u|^{\frac{2}{1-\beta}}}-1\right)\mathrm{d}x=+\infty.
\end{equation*}
This ends the proof of Theorem \ref{conclu1}.

\section{Proof of Theorem \ref{conclu2}}
In this section, we prove the improvement of the Adams' inequality using the inequality we have proved in the previous section, that is the proof of Theorem \ref{conclu2}.

Let $\{u_n\}\subset X_\beta$ and $u\in X_\beta\setminus\{0\}$ be such that
${\|u_n\|}_\beta=1$ and $u_n\rightharpoonup u$ weakly in $X_\beta$.

\noindent For $a,b\in \mathbb{R}$, $q>1$. If $q'$ its conjugate i.e. $\frac{1}{q}+\frac{1}{q'}=1$ and by Young inequality, we have
\begin{equation*}
  e^{a+b}\leq \frac{1}{q}e^{qa}+\frac{1}{q'}e^{q'b}.
\end{equation*}
Also we have
\begin{equation*}
  \begin{split}
    |u_n|^{\frac{2}{1-\beta}}=|u_n-u+u|^{\frac{2}{1-\beta}}
    & \leq\big(|u_n-u|+|u|\big)^{\frac{2}{1-\beta}} \\
    &\leq(1+\epsilon)|u_n-u|^{\frac{2}{1-\beta}}+\left(1-\frac{1}{(1+\epsilon)^{\frac{1-\beta}{1+\beta}}}\right)
    ^{\frac{\beta+1}{\beta-1}}|u|^{\frac{2}{1-\beta}},
  \end{split}
\end{equation*}
which implies that
\begin{align}\label{3.1}
  \int^{}_{\mathbb{R}^4}\left(e^{p\alpha_\beta|u_n|^{\frac{2}{1-\beta}}}-1\right)\mathrm{d}x
  & \leq\frac{1}{q}\int^{}_{\mathbb{R}^4}\left(e^{(1+\epsilon)pq\alpha_\beta|u_n-u|^{\frac{2}{1-\beta}}}-1\right)\mathrm{d}x \notag\\
  & +\frac{1}{q'}\int^{}_{\mathbb{R}^4}\left(e^{pq'\alpha_\beta
  \left(1-\frac{1}{(1+\epsilon)^{\frac{1-\beta}{1+\beta}}}\right)^{\frac{\beta+1}{\beta-1}}|u|^{\frac{2}{1-\beta}}}-1\right)\mathrm{d}x.
\end{align}
By (\ref{therom1.2(1)}), we have
\begin{align}\label{3.2}
   & \int^{}_{\mathbb{R}^4}\left(e^{pq'\alpha_\beta\left(1-\frac{1}{(1+\epsilon)^{\frac{1-\beta}{1+\beta}}}\right)^{\frac{\beta+1}{\beta-1}}
    |u|^{\frac{2}{1-\beta}}}-1\right)\mathrm{d}x  \notag\\
   &\leq\int^{}_{\mathbb{R}^4}\left(e^{\left|\left(pq'\alpha_\beta\left(1-\frac{1}{(1+\epsilon)^{\frac{1-\beta}{1+\beta}}}\right)^{\frac{\beta+1}{\beta-1}}\right)
    ^{\frac{1-\beta}{2}}|u|\right|^{\frac{2}{1-\beta}}}-1\right)\mathrm{d}x<+\infty.
\end{align}
To complete the proof, we have to prove that for every $p$ such that $0<p<P_\beta(u)$,
\begin{equation}\label{3.3}
  \sup_n\int^{}_{\mathbb{R}^4}\left(e^{(1+\epsilon)pq\alpha_\beta|u_n-u|^{\frac{2}{1-\beta}}}-1\right)\mathrm{d}x<+\infty,
\end{equation}
for some $q>1$ and $\epsilon>0$.

If $\|u\|_\beta<1$ and $0<p<P_\beta(u)=\frac{1}{\left(1-\|u\|^2_\beta\right)^{\frac{1}{1-\beta}}}$, there exists $\nu>0$ such that
\begin{equation*}
  p\left(1-\|u\|^2_\beta\right)^{\frac{1}{1-\beta}}(1+\nu)<1.
\end{equation*}
Moreover, by Brezis-Lieb Lemma \cite{Brezis}, we can easily get
\begin{equation*}
  \|u_n-u\|^2_\beta=\|u_n\|^2_\beta-\|u\|^2_\beta+o(1)=1-\|u\|^2_\beta+o(1),
\end{equation*}
where $o(1)\rightarrow0$ as $n\rightarrow+\infty$.

\noindent Then
\begin{equation*}
  \lim_{n\rightarrow+\infty}\|u_n-u\|^{\frac{2}{1-\beta}}_\beta=\left(\lim_{n\rightarrow+\infty}\|u_n-u\|^2_\beta\right)^{\frac{1}{1-\beta}}
  =\left(1-\|u\|^2_\beta\right)^{\frac{1}{1-\beta}}.
\end{equation*}
So $\text{for all} \,\, \epsilon >0$, there exists $N_\epsilon\geq1$ such that
\begin{equation*}
  \|u_n-u\|^{\frac{2}{1-\beta}}_\beta \leq\left(1-\|u\|^2_\beta\right)^{\frac{1}{1-\beta}}(1+\epsilon),\quad\text{for all}  \,\,\, n\geq N_\epsilon.
\end{equation*}
Choose $q=1+\epsilon$, $\epsilon=\sqrt[3]{1+\nu}-1$, one has
\begin{equation*}
  pq(1+\epsilon)\|u_n-u\|^{\frac{2}{1-\beta}}_\beta \leq p(1+\nu)\left(1-\|u\|^2_\beta\right)^{\frac{1}{1-\beta}}<1.
\end{equation*}
Consequently, using (\ref{therom1.2(2)}), we get
\begin{equation}\label{3.4}
  \begin{split}
    \int^{}_{\mathbb{R}^4}\left(e^{(1+\epsilon)pq\alpha_\beta|u_n-u|^{\frac{2}{1-\beta}}}-1\right)\mathrm{d}x
    & =\int^{}_{\mathbb{R}^4}\left(e^{pq(1+\epsilon)\|u_n-u\|^{\frac{2}{1-\beta}}_\beta\alpha_\beta
    \left(\frac{|u_n-u|}{\|u_n-u\|_\beta}\right)^{\frac{2}{1-\beta}}}-1\right)\mathrm{d}x \\
    & <+\infty.
  \end{split}
\end{equation}

If $\|u\|_\beta=1$ and $0<p<P_\beta(u)=+\infty$, we also have
\begin{equation*}
  \lim_{n\rightarrow+\infty}\|u_n-u\|^{\frac{2}{1-\beta}}_\beta=\left(\lim_{n\rightarrow+\infty}\|u_n-u\|^2_\beta\right)^{\frac{1}{1-\beta}}=0.
\end{equation*}
So $\text{for all}  \,\, \epsilon >0$, there exists $N_\epsilon\geq1$ such that
\begin{equation*}
  \|u_n-u\|^{\frac{2}{1-\beta}}_\beta<\frac{1}{p(1+\epsilon)^2},\quad\text{for all}  \,\,\, n\geq N_\epsilon.
\end{equation*}
Choose $q=1+\epsilon$, one has
\begin{equation*}
  pq(1+\epsilon)\|u_n-u\|^{\frac{2}{1-\beta}}_\beta < p(1+\epsilon)^2\cdot\frac{1}{p(1+\epsilon)^2}=1.
\end{equation*}
Consequently, using (\ref{therom1.2(2)}), we also get
\begin{equation}\label{3.5}
  \begin{split}
    \int^{}_{\mathbb{R}^4}\left(e^{(1+\epsilon)pq\alpha_\beta|u_n-u|^{\frac{2}{1-\beta}}}-1\right)\mathrm{d}x
    & =\int^{}_{\mathbb{R}^4}\left(e^{pq(1+\epsilon)\|u_n-u\|^{\frac{2}{1-\beta}}_\beta\alpha_\beta
    \left(\frac{|u_n-u|}{\|u_n-u\|_\beta}\right)^{\frac{2}{1-\beta}}}-1\right)\mathrm{d}x \\
    & <+\infty.
  \end{split}
\end{equation}
Now (\ref{3.3}) follows from (\ref{3.4}) and (\ref{3.5}).

Combining (\ref{3.1}), (\ref{3.2}) and (\ref{3.3}), the proof is complete.

\section{Application}

In this section, we introduce an application of the previous inequalities and give the details of the proof of Theorem \ref{conclu3}.

Notice that problem (\ref{problem}) has a variational structure, so we use the variational methods to deal with the problem. The energy function $I_{\beta}$ associated with problem (\ref{problem}) is given by
\begin{equation*}
  I_{\beta}(v)=\frac{1}{2}\| v\|^2_{\beta}-\int^{}_{\mathbb{R}^4}F(v)dx,\quad\text{for all}  \,\,\, v\in X_\beta.
\end{equation*}
We can easily see that $I_\beta$ is well defined and $I_\beta\in C^1(X_\beta,\mathbb{R})$.

A weak solution of problem (\ref{problem}) is a function $u\in X_\beta$ which satisfies
\begin{equation*}
  \int_{\mathbb{R}^4}\omega_\beta(x)\triangle u\triangle v\mathrm{d}x=\int_{\mathbb{R}^4}f(u)v\mathrm{d}x, \quad \text{for all}  \,\,\, v\in X_\beta.
\end{equation*}
It is well-known that finding a weak solution of problem (\ref{problem}) is equivalent to finding a critical point of $I_\beta$.

First, we give the embedding properties which will be useful in the rest of our work.
\begin{proposition}\label{embedding}
Let $0<\beta<1$. Assume that $(D_1)$ and $(D_2)$ holds. Then for all $q\geq\frac{2}{1-\beta}$, $X_\beta\hookrightarrow L^q(\mathbb{R}^4)$ with continuous embedding. Furthermore, for all $q>\frac{2}{1-\beta}$, the embedding $X_\beta\hookrightarrow L^q(\mathbb{R}^4)$ is compact.
\end{proposition}
\begin{proof}
Let $u\in X_\beta$, by Lemma \ref{radial}, we have
\begin{equation*}
  |u(e_1)|\leq \frac{\sqrt{2}c_1}{4\pi}\|u\|_\beta\left(\int^{+\infty}_{1}\chi(y)^{-1}y^{-1}\mathrm{d}y\right)^{\frac{1}{2}}+c_2.
\end{equation*}
Set
\begin{equation*}
  v(x)=\begin{cases}
  u(x)-u(e_1),&0\leq|x|<1,\\
  0,&|x|\geq1,
  \end{cases}
\end{equation*}
where $e_1=(1,0,0,0)\in \mathbb{R}^4$. Clearly $v\in W^{2,2}_{0,rad}(B,\omega)$.

\noindent From \cite[Lemma 2.1]{ZW}, we have
\begin{equation*}
  |v(x)|\leq\frac{1}{2\sqrt{2}\pi}\left(\frac{\left(\log{\frac{e}{|x|}}\right)^{1-\beta}-1}{1-\beta}\right)^{\frac{1}{2}}\|u\|_\beta.
\end{equation*}
Then
\begin{equation*}
  \begin{split}
    |u(x)| & \leq|v(x)|+|u(e_1)| \\
    & \leq\frac{1}{2\sqrt{2}\pi}\left(\frac{\left(\log{\frac{e}{|x|}}\right)^{1-\beta}-1}{1-\beta}\right)^{\frac{1}{2}}\|u\|_\beta
    +\frac{\sqrt{2}c_1}{4\pi}\|u\|_\beta\left(\int^{+\infty}_{1}\chi(y)^{-1}y^{-1}\mathrm{d}y\right)^{\frac{1}{2}}+c_2, \quad \text{for} \,\, |x|<1.
  \end{split}
\end{equation*}
By (\ref{chi0}), it is easy to infer that
\begin{equation}\label{pro2.1}
  \int^{}_{|x|<1}|u(x)|^{\frac{2k}{1-\beta}}\mathrm{d}x \leq c_3\|u\|^{\frac{2k}{1-\beta}}_\beta+c_4 \leq(c_3+c_5)\|u\|^{\frac{2k}{1-\beta}}_\beta.
\end{equation}
On the other hand, we can infer from (\ref{2.3}) that
\begin{equation}\label{pro2.2}
  \int^{}_{|x|\geq1}|u|^{\frac{2k}{1-\beta}}\mathrm{d}x \leq c_6\|u\|^{\frac{2k}{1-\beta}}_\beta+c_7 \leq(c_6+c_8)\|u\|^{\frac{2k}{1-\beta}}_\beta.
\end{equation}
Combining (\ref{pro2.1}) and (\ref{pro2.2}), one has
\begin{equation*}
  \int^{}_{\mathbb{R}^4}|u|^{\frac{2k}{1-\beta}}\mathrm{d}x \leq c_9\|u\|_\beta^{\frac{2k}{1-\beta}}.
\end{equation*}
Then we deduce that the continuous embedding $X_\beta\hookrightarrow L^{\frac{2k}{1-\beta}}(\mathbb{R}^4)$ holds for all $k\in N^*$.
 By simple interpolation, it follows that $X_\beta$ is continuously embedded into $L^q(\mathbb{R}^4)$,$\, \text{for all} \,\, q\geq\frac{2}{1-\beta}$.
 Proceeding as the proof of \cite[proposition 2]{SRn}, we deduce that the embedding $X_\beta\hookrightarrow L^q(\mathbb{R}^4)$ is compact for all $q>\frac{2}{1-\beta}$.
\end{proof}

Now we consider problem (\ref{problem}) and verify that the corresponding energy function $I_\beta$ has the Mountain Pass geometry structure.
\begin{lemma}\label{mp}
Assume that $(F_1)-(F_3)$ hold. Then

(i) There exist $\rho_0,\tau>0$ such that $ I_{\beta}(u)\geq \tau$, $\text{for} \,\, u\in X_{\beta}$, $\|u\|_{\beta}=\rho_0$.

(ii) There exists $\|e_0\|_{\beta}>\rho_0$ such that $ I_{\beta}(e_0)<0$.
\end{lemma}
\begin{proof}
(i) Given $u\in X_\beta$, by $(F_1)$, we get
\begin{equation}\label{Lem2.1}
  \begin{split}
     \int^{}_{\mathbb{R}^4}|F(u)|\mathrm{d}x
     & \leq c_0\left(\int^{}_{\mathbb{R}^4}|u|^p\mathrm{d}x +\int^{}_{\mathbb{R}^4}|u|^q\left(e^{\alpha |u|^{\frac{2}{1-\beta}}}-1\right)\mathrm{d}x \right) \\
     & \leq c_0\left(|u|^p_{L^p(\mathbb{R}^4)}+|u|^q_{L^{q+1}(\mathbb{R}^4)} \left(\int^{}_{\mathbb{R}^4}\left(e^{\alpha(q+1)|u|^{\frac{2}{1-\beta}}}-1\right)
     \mathrm{d}x \right)^{\frac{1}{q+1}} \right).
  \end{split}
\end{equation}
For $u\in X_\beta$ such that
\begin{equation*}
  \|u\|_\beta=\rho_0< \left(\frac{\alpha_\beta}{\alpha(1+q)}\right)^{\frac{1-\beta}{2}}.
\end{equation*}
By (\ref{therom1.2(2)}), there exists a constant $C_3>0$ such that
\begin{equation*}
  \int^{}_{\mathbb{R}^4}\left( e^{\alpha(q+1)|u|^{\frac{2}{1-\beta}}}-1 \right)\mathrm{d}x
  =\int^{}_{\mathbb{R}^4}\left( e^{\alpha(q+1) \|u\|^{\frac{2}{1-\beta}}_{\beta} \left|\frac{u}{\|u\|_{\beta}}\right|^{\frac{2}{1-\beta}}}-1 \right)\mathrm{d}x \leq C_3.
\end{equation*}
Since $X_\beta\hookrightarrow L^p\left(\mathbb{R}^4\right)$, $X_\beta\hookrightarrow L^{q+1}\left(\mathbb{R}^4\right)$ and combine with (\ref{Lem2.1}), we have
\begin{equation*}
  \int^{}_{\mathbb{R}^4}\left|F(u)\right|\mathrm{d}x \leq C_4\left(\|u\|^p_\beta+\|u\|^q_\beta\right)
\end{equation*}
where $C_4>0$. Then
\begin{equation*}
  I_\beta(u)=\frac{1}{2}\|u\|^2_{\beta}-\int^{}_{\mathbb{R}^4}F(u)\mathrm{d}x \geq \frac{1}{2}\|u\|^2_{\beta}
  -C_4\left(\|u\|^p_{\beta}+\|u\|^q_{\beta}\right),\quad\text{for}  \,\,\, u\in X_\beta,\, \|u\|_\beta=\rho_0.
\end{equation*}
One can choose $\rho_0>0$ small enough such that
\begin{equation*}
  \rho_0<\min\left\{\left(\frac{\alpha_\beta}{\alpha(1+q)}\right)^{\frac{1-\beta}{2}},\left(\frac{1}{4C_4}\right)^{\frac{1}{\min\{p,q\}-2}}\right\}.
\end{equation*}
Then we deduce that
\begin{equation*}
  I_\beta(u) \geq \rho^2_0 \left(\frac{1}{2}-2C_4\rho^{\min\{p,q\}-2}_0\right)=\tau>0,\quad\text{for}  \,\,\, u\in X_\beta,\, \|u\|_\beta=\rho_0.
\end{equation*}
(ii)From $(F_3)$, we have
\begin{equation*}
  \lim_{s\rightarrow\infty}\frac{f(s)s}{e^{\alpha s^{\frac{2}{1-\beta}}}}=\lim_{s\rightarrow\infty}\frac{f(s)s}{s^{\frac{2}{1-\beta}}}
  \cdot\frac{s^{\frac{2}{1-\beta}}}{e^{\alpha s^{\frac{2}{1-\beta}}}}\geq\kappa>0.
\end{equation*}
Since $\lim\limits_{s\rightarrow\infty}\frac{s^{\frac{2}{1-\beta}}}{e^{\alpha s^{\frac{2}{1-\beta}}}}=0$, we infer
\begin{equation*}
  \lim_{s\rightarrow\infty}\frac{f(s)}{s^{\frac{2}{1-\beta}-1}}=\infty.
\end{equation*}
Then, there exist $s_1>0$ and $M>0$ such that for any $s\geq s_1$,
\begin{equation*}
  f(s)\geq Ms^{\frac{2}{1-\beta}-1}.
\end{equation*}
It follows that
\begin{equation*}
  F(s)\geq \frac{1-\beta}{2}Ms^{\frac{2}{1-\beta}}.
\end{equation*}
By $(F_1)$, we have
\begin{equation*}
 \begin{split}
    I_\beta(tu) & =\frac{1}{2}t^2\|u\|^2_{\beta}-\int^{}_{\mathbb{R}^4}F(tu)\mathrm{d}x \\
      & \leq \frac{1}{2}t^2\|u\|^2_{\beta} -\int_{|tu|\geq s_1}F(tu)\mathrm{d}x \\
      & \leq \frac{1}{2}t^2\|u\|^2_{\beta} -\frac{1-\beta}{2}Mt^{\frac{2}{1-\beta}}\int_{|tu|\geq s_1}u^{\frac{2}{1-\beta}}\mathrm{d}x.
 \end{split}
\end{equation*}
Then, we obtain
\begin{equation*}
  I_{\beta}(tu)\rightarrow -\infty,\quad \text{when}\,\, t\rightarrow +\infty.
\end{equation*}
Thus, there exists $e_0\in X_{\beta}$, $\|e_0\|_{\beta}>\rho_0$ such that $ I_{\beta}(e_0)<0$.
\end{proof}

By Lemma \ref{mp} and the Mountain Pass Theorem, we know that there exists a $(PS)_m$ sequence of the functional $I_\beta$ at the mountain pass level
\begin{equation*}
  m=\inf\limits_{\gamma\in\Gamma}\sup\limits_{0\leq t\leq1} I_{\beta}\big(\gamma(t)\big),
\end{equation*}
where
\begin{equation*}
  \Gamma=\left\{\gamma\in C\left([0,1],X_{\beta}\right),\gamma(0)=0,I_{\beta}\big(\gamma(1)\big)<0\right\}.
\end{equation*}

Next, we establish an important estimate involving the level $m$.
\begin{lemma}
Assume that $(F_1)-(F_3)$ hold. Then there holds
\begin{equation*}
  0<m<m^{*}\triangleq\frac{1}{2}\left(\frac{\alpha_{\beta}}{\alpha}\right)^{1-\beta}.
\end{equation*}
\end{lemma}
\begin{proof}
By the previous lemma, clearly, $m\geq\rho_1>0$.

Next, we consider the family of Adams type function
\begin{equation*}
  \omega_n(x)=\begin{cases}
   \left(\frac{\log (e^4n)}{\alpha_\beta}\right)^{\frac{1-\beta}{2}}
   -\frac{|x|^{2(1-\beta)}}{2\left(\frac{\alpha_\beta}{4n}\right)^{\frac{1-\beta}{2}}\left(\frac{1}{4}\log(e^4n)\right)^{\frac{1+\beta}{2}}}
   +\frac{1}{2\left(\frac{\alpha_\beta}{4}\right)^{\frac{1-\beta}{2}}\left(\frac{1}{4}\log(e^4n)\right)^{\frac{1+\beta}{2}}}
   ,&0\leq|x|\leq \frac{1}{\sqrt[4]{n}},\\
   \frac{\left(\log\frac{e}{|x|}\right)^{1-\beta}} {\left(\frac{\alpha_\beta}{16}\log(e^4n)\right)^{\frac{1-\beta}{2}}}
   ,&\frac{1}{\sqrt[4]{n}}<|x|\leq\frac{1}{2},\\
   \eta_n,&|x|>\frac{1}{2},
  \end{cases}
\end{equation*}
where $\eta_n\in C^\infty_0(B_1)$ is such that
\begin{equation*}
  \eta_n|_{x=\frac{1}{2}}=\frac{1}{\left(\frac{\alpha_\beta}{16}\log(e^4n)\right)^{\frac{1-\beta}{2}}}(\log2e)^{1-\beta},
  \frac{\partial \eta_n}{\partial x}|_{x=\frac{1}{2}}
  =-\frac{2(1-\beta)}{\left(\frac{\alpha_\beta}{16}\log(e^4n)\right)^{\frac{1-\beta}{2}}}(\log2e)^{-\beta},
\end{equation*}
and $\eta_n,\nabla\eta_n,\triangle\eta_n$ are all $O\left(\frac{1}{\left(\log(e^4n)\right)^\frac{1-\beta}{2}}\right)$ as $n\rightarrow\infty$.

\noindent As in the proof of Section 3, after a series of computation, we get
\begin{equation*}
  \|\omega_n\|_\beta\rightarrow1,\quad \text{as} \,\, n\rightarrow\infty.
\end{equation*}
Let
\begin{equation*}
  v_n=\frac{\omega_n}{\|\omega_n\|_\beta}.
\end{equation*}
We also have
\begin{equation*}
  \|v_n\|_\beta=1,
\end{equation*}
and
\begin{equation}\label{psc1}
  v_n(x)\geq\left(\frac{\log(e^4n)}{\alpha_\beta}\right)^{\frac{1-\beta}{2}}+o(1),\quad\text{for all}  \,\,\, x\in\mathbb{R}^4, \,\, 0\leq|x|\leq\frac{1}{\sqrt[4]{n}}.
\end{equation}

In order to prove $m<m^{*}$, we only need to prove that
\begin{equation*}
  \max\limits_{t\geq0}I_\beta(tv_n)<m^*.
\end{equation*}
By contradiction, we suppose that
\begin{equation*}
  \max\limits_{t\geq0}I_\beta(tv_n)\geq m^*.
\end{equation*}
Then, there exists $t_n>0$ such that
\begin{equation*}
  \max\limits_{t\geq0}I_\beta(tv_n)=I_\beta(t_nv_n)\geq m^*,
\end{equation*}
that is
\begin{equation*}
  \frac{1}{2}t_n^2\|v_n\|^2_{\beta}-\int_{\mathbb{R}^4}F(t_nv_n)\mathrm{d}x\geq\frac{1}{2}\left(\frac{\alpha_{\beta}}{\alpha}\right)^{1-\beta}.
\end{equation*}
Combining with $(F_1)$, we get
\begin{equation}\label{psc2}
  t_n^2\geq\left(\frac{\alpha_{\beta}}{\alpha}\right)^{1-\beta}.
\end{equation}
Moreover, we have
\begin{equation*}
  \frac{d}{dt}I_\beta(tv_n)\big|_{t=t_n}=t_n-\int_{\mathbb{R}^4}f(t_nv_n)v_n\mathrm{d}x=0,
\end{equation*}
that is
\begin{equation}\label{psc3}
  t^2_n=\int_{\mathbb{R}^4}f(t_nv_n)t_nv_n\mathrm{d}x.
\end{equation}
Now we claim that $\{t_n\}$ is bounded.

\noindent In fact, it follows from $(F_3)$ that for any $\epsilon>0$, there exists $s(\epsilon)$ such that
\begin{equation}\label{psc4}
  f(s)s\geq(\kappa-\epsilon)e^{\alpha s^{\frac{2}{1-\beta}}}, \quad \text{for all}  \,\,\, s\geq s(\epsilon).
\end{equation}
Notice that $t_nv_n\geq s(\epsilon)$ for $0\leq|x|\leq \frac{1}{\sqrt[4]{n}}$ and for large $n\in \mathbb{N}$. Then it follows from (\ref{psc1}), (\ref{psc3}) and (\ref{psc4}) that for large $n$,
\begin{equation}\label{psc5}
  \begin{split}
    t^2_n & =\int_{\mathbb{R}^4}f(t_nv_n)t_nv_n\mathrm{d}x \\
          & \geq \int_{0\leq|x|\leq \frac{1}{\sqrt[4]{n}}}f(t_nv_n)t_nv_n\mathrm{d}x \\
          & \geq (\kappa-\epsilon)\int_{0\leq|x|\leq \frac{1}{\sqrt[4]{n}}} e^{\alpha t_n^{\frac{2}{1-\beta}}v_n^{\frac{2}{1-\beta}}}\mathrm{d}x \\
          & \geq 2\pi^2(\kappa-\epsilon)\int^{\frac{1}{\sqrt[4]{n}}}_{0} r^3e^{\alpha t_n^{\frac{2}{1-\beta}}\frac{\log(e^4n)}{\alpha_\beta}}\mathrm{d}r,
  \end{split}
\end{equation}
that is
\begin{equation*}
  1\geq2\pi^2(\kappa-\epsilon)e^{\alpha\frac{\log(e^4n)}{\alpha_\beta}t_n^{\frac{2}{1-\beta}}-\log(4n)-\log t_n^2}.
\end{equation*}
Therefore, $\{t_n\}$ is bounded.

\noindent From (\ref{psc2}), we have
\begin{equation*}
  \lim_{n\rightarrow+\infty}t_n^2\geq\left(\frac{\alpha_{\beta}}{\alpha}\right)^{1-\beta}.
\end{equation*}
In fact, we can prove that
\begin{equation}\label{psc6}
  \lim_{n\rightarrow+\infty}t_n^2=\left(\frac{\alpha_{\beta}}{\alpha}\right)^{1-\beta}.
\end{equation}
Indeed, if $\lim\limits_{n\rightarrow+\infty}t_n^2>\left(\frac{\alpha_{\beta}}{\alpha}\right)^{1-\beta}$, it follows that for large $n$ there exists $\delta>0$ such that
\begin{equation}\label{psc7}
  t_n^{\frac{2}{1-\beta}}\geq\frac{\alpha_\beta}{\alpha}+\delta.
\end{equation}
Combining (\ref{psc5}) and (\ref{psc7}), we infer $t_n\rightarrow\infty$, as $n\rightarrow\infty$ which contradicts the boundedness of $\{t_n\}$.

\noindent Therefore, (\ref{psc6}) follows immediately.

\noindent Combining (\ref{psc5}) and (\ref{psc6}), we have
\begin{equation*}
  \begin{split}
    \lim\limits_{n\rightarrow\infty}t_n^2=\left(\frac{\alpha_{\beta}}{\alpha}\right)^{1-\beta}
    & \geq(\kappa-\epsilon)\lim\limits_{n\rightarrow\infty}\int_{0\leq|x|\leq\frac{1}{\sqrt[4]{n}}}
      e^{\alpha t_n^{\frac{2}{1-\beta}}v_n^{\frac{2}{1-\beta}}}\mathrm{d}x \\
    & \geq2\pi^2(\kappa-\epsilon)\lim\limits_{n\rightarrow\infty}\int^{\frac{1}{\sqrt[4]{n}}}_{0}r^3 e^{\log(e^4n)}\mathrm{d}r \\
    & =\frac{\pi^2}{2}(\kappa-\epsilon)e^4.
  \end{split}
\end{equation*}
It follows that
\begin{equation*}
  \kappa\leq\frac{2}{\pi^2e^4}\left(\frac{\alpha_\beta}{\alpha}\right)^{1-\beta},
\end{equation*}
which contradicts with $(F_3)$.
\end{proof}

\begin{lemma}\label{bounded}
Assume that $(F_1)-(F_3)$ hold. Let $\{u_n\}$ be a $(PS)_m$ sequence of $I_\beta$ with $m\in(0,m^{*})$. Then $\{u_n\}$ is bounded in $X_\beta$.
\end{lemma}
\begin{proof}
Let $\{u_n\}\subset X_\beta$ be a $(PS)_m$ sequence of $I_\beta$ with $m\in(0,m^{*})$. It follows that
\begin{equation}\label{I}
  I_{\beta}(u_n)=\frac{1}{2}\| u_n\|^2_{\beta}-\int^{}_{\mathbb{R}^4}F(u_n)dx\rightarrow m,\quad n\rightarrow+\infty,
\end{equation}
and
\begin{equation}\label{I1}
  \left|I'_{\beta}(u_n)u_n\right|=\left|\|u_n\|^2_\beta-\int_{\mathbb{R}^4}f(u_n)u_n\mathrm{d}x\right|\leq \epsilon_n\|u_n\|_\beta,
\end{equation}
where $\epsilon_n\rightarrow0$, as $n\rightarrow+\infty$.

\noindent From $(F_2)$, we have for any $\epsilon>0$ there exists $s_\epsilon>0$ such that
\begin{equation}\label{Fs}
  F(s)\leq\epsilon sf(s), \quad\text{for all}  \,\,\, s\geq s_\epsilon.
\end{equation}
Combining (\ref{I}), (\ref{I1}) and (\ref{Fs}), for fix $0<\epsilon<\frac{1}{2}$, we have
\begin{equation*}
  \begin{split}
     \frac{1}{2}\|u_n\|^2_\beta
     & \leq C+\int^{}_{\mathbb{R}^4}F(u_n)dx \\
     & \leq C+\int^{}_{|u_n|\leq s_\epsilon}F(u_n)dx+\int^{}_{|u_n|\geq s_\epsilon}F(u_n)dx \\
     & \leq C'+\epsilon\int_{\mathbb{R}^4}f(u_n)u_n\mathrm{d}x \\
     & \leq C'+\epsilon\|u_n\|^2_\beta+\epsilon\epsilon_n\|u_n\|_\beta.
   \end{split}
\end{equation*}
Therefore, we can deduce that $\{u_n\}$ is bounded in $X_\beta$.
\end{proof}

Now we turn to the compactness properties of $(PS)_m$ sequence.
\begin{lemma}\label{convergence1}
Assume that $(F_1)-(F_3)$ hold. Let $\{u_n\}$ be a $(PS)_m$ sequence of $I_\beta$ such that
\begin{equation}\label{Lemma5}
  \liminf_{n\rightarrow+\infty}\|u_n\|_{\beta} <\left(\frac{\alpha_{\beta}}{\alpha}\right)^{\frac{1-\beta}{2}},
\end{equation}
then $\{u_n\}$ has a convergence subsequence in $X_\beta$.
\end{lemma}
\begin{proof}
By Lemma \ref{bounded}, up to a subsequence, one has
\begin{equation*}
  u_n\rightharpoonup u \;\; \text{weakly} \; \text{in} \; X_\beta \;\; \text{and} \;\; u_n(x)\rightarrow u(x) \; a.e. \; x\in \mathbb{R}^4.
\end{equation*}

\noindent By (\ref{Lemma5}), up to a subsequence, there exist $r_1>1$ and $0<\gamma_0<1$ such that
\begin{equation}\label{Lemma5(1)}
  \|u_n\|_{\beta} < \left(\frac{\gamma_0\alpha_{\beta}}{r_1\alpha}\right)^{\frac{1-\beta}{2}}.
\end{equation}
By $(F_1)$, we get
\begin{equation}\label{Lemma5(2)}
  \left|\int^{}_{\mathbb{R}^4}f(u_n)(u_n-u)\mathrm{d}x\right|
  \leq c_1\int^{}_{\mathbb{R}^4}|u_n|^{p-1}|u_n-u|\mathrm{d}x +c_1\int^{}_{\mathbb{R}^4}|u_n|^{q-1}\left(e^{\alpha|u_n|^{\frac{2}{1-\beta}}}-1\right)|u_n-u|\mathrm{d}x.
\end{equation}
Using H\"{o}lder inequality, we have
\begin{equation}\label{Lemma5(3)}
 \begin{split}
   \int^{}_{\mathbb{R}^4}|u_n|^{p-1}|u_n-u|\mathrm{d}x
   & \leq \left( \int^{}_{\mathbb{R}^4}|u_n|^{p}\right)^{\frac{p-1}{p}} \left( \int^{}_{\mathbb{R}^4}|u_n-u|^{p}\right)^{\frac{1}{p}} \\
   & =|u_n|^{p-1}_{L^{p}(\mathbb{R}^4)} |u_n-u|_{L^{p}(\mathbb{R}^4)}.
 \end{split}
\end{equation}
Let $r_2>1$ and $r_3>1$ be such that

\begin{align*}
\left\{
  \begin{aligned}
     & \frac{1}{r_1}+\frac{1}{r_2}+\frac{1}{r_3}=1, \\
     & r_2>\frac{r_1}{r_1-1}, \\
     & r_2\geq \frac{2}{(1-\beta)(q-1)}, \\
     & r_3>\frac{2}{1-\beta}.
  \end{aligned}
\right.
\end{align*}
From (\ref{therom1.2(2)}), (\ref{Lemma5(1)}) and H\"{o}lder inequality, we have
\begin{equation}\label{Lemma5(4)}
  \begin{split}
   &\int^{}_{\mathbb{R}^4}|u_n|^{q-1} \left(e^{\alpha|u_n|^{\frac{2}{1-\beta}}}-1\right)|u_n-u|\mathrm{d}x \\
   & \leq \left|u_n\right|^{q-1}_{L^{r_2(q-1)}(\mathbb{R}^4)} |u_n-u|_{L^{r_3}(\mathbb{R}^4)}
     \left(\int^{}_{\mathbb{R}^4}\left(e^{r_1\alpha|u_n|^{\frac{2}{1-\beta}}}-1\right)\mathrm{d}x\right)^{\frac{1}{r_1}} \\
   & \leq |u_n|^{q-1}_{L^{r_2(q-1)}(\mathbb{R}^4)} |u_n-u|_{L^{r_3}(\mathbb{R}^4)}
     \left(\int^{}_{\mathbb{R}^4}\left(e^{r_1\alpha\|u_n\|_{\beta}^{\frac{2}{1-\beta}} \left(\frac{u_n}{\|u_n\|_{\beta}}\right)^{\frac{2}{1-\beta}}}-1\right)\mathrm{d}x\right)^{\frac{1}{r_1}} \\
   & \leq |u_n|^{q-1}_{L^{r_2(q-1)}(\mathbb{R}^4)} |u_n-u|_{L^{r_3}(\mathbb{R}^4)}
     \left(\int^{}_{\mathbb{R}^4}\left(e^{\gamma_0\alpha_{\beta}\left(\frac{u_n}{\|u_n\|_{\beta}}\right)^{\frac{2}{1-\beta}}}-1\right)\mathrm{d}x\right)^{\frac{1}{r_1}} \\
   & \leq c_2|u_n|^{q-1}_{L^{r_2(q-1)}(\mathbb{R}^4)} |u_n-u|_{L^{r_3}(\mathbb{R}^4)}.
  \end{split}
\end{equation}
Since the embeddings $X_\beta\hookrightarrow L^{p}\left(\mathbb{R}^4\right)$ and $X_\beta\hookrightarrow L^{r_3}\left(\mathbb{R}^4\right)$ are compact, we can immediately deduce from (\ref{Lemma5(2)}), (\ref{Lemma5(3)}) and (\ref{Lemma5(4)}) that
\begin{equation*}
  \int^{}_{\mathbb{R}^4}f(u_n)(u_n-u)\mathrm{d}x\rightarrow0.
\end{equation*}
On the other hand, since $\left\langle I'_{\beta}(u_n),u_n-u\right\rangle\rightarrow0$, as $n\rightarrow+\infty$, then we infer
\begin{equation*}
  \int^{}_{\mathbb{R}^4}\omega_{\beta}\triangle u_n\triangle(u_n-u)\rightarrow0,\quad n\rightarrow+\infty.
\end{equation*}
Thus,
\begin{equation*}
  \|u_n-u\|^2_\beta =\langle u_n,u_n-u\rangle -\langle u,u_n-u\rangle =\langle u_n,u_n-u\rangle+o_n(1)\rightarrow0,\quad n\rightarrow+\infty.
\end{equation*}
That is, $u_n\rightarrow u$ strongly in $X_\beta$.
\end{proof}

\begin{lemma}\label{end}
Assume that $(F_1)-(F_3)$ hold. Then $I_\beta$ satisfies the $(PS)_m$ condition if $m\in(0,m^{*})$.
\end{lemma}
\begin{proof}
Let $\{u_n\}\subset X_\beta$ be a $(PS)_m$ sequence of $I_\beta$ with $m\in\left(0,m^{*}\right)$. By Lemma \ref{bounded}, up to a subsequence, one has
\begin{equation*}
  u_n\rightharpoonup u \;\; \text{weakly} \; \text{in} \; E_\beta \;\; \text{and} \;\; u_n(x)\rightarrow u(x)\quad a.e. \,\, x\in \mathbb{R}^4,
\end{equation*}
and
\begin{equation*}
  \left|\left\langle I_{\beta}'(u_n),u_n\right\rangle\right|=\left|\|u_n\|^2_{\beta}-\int^{}_{\mathbb{R}^4}f(u_n)u_n \mathrm{d}x\right| \leq C,\quad\text{for all}\,\, n.
\end{equation*}
Since $u_n$ is bounded in $X_\beta$, then
\begin{equation*}
  \int^{}_{\mathbb{R}^4}f(u_n)u_n \mathrm{d}x \leq C,\quad\text{for all}  \,\, n.
\end{equation*}
By \cite[Lemma 2.1]{FMB}, we get $f(u_n)\rightarrow f(u)$ in $L^1_{loc}(\mathbb{R}^4)$. Then, it follows from $(F_2)$ and Generalized Lebesgue Dominated Convergence Theorem that  $F(u_n)\rightarrow F(u)$ in $L^1_{loc}(\mathbb{R}^4)$. Moreover, inspired by \cite{Fun}, we can deduce that
\begin{equation*}
  F(u_n)\rightarrow F(u) \quad \text{in} \,\, L^1(\mathbb{R}^4) \quad \text{as} \,\, n\rightarrow+\infty.
\end{equation*}
In fact, by (\ref{Fs}), we have
\begin{equation}\label{Fun1}
  \int_{|u_n|\geq s_\epsilon}|F(u_n)|\mathrm{d}x\leq\epsilon \int_{|u_n|\geq s_\epsilon}f(u_n)u_n\mathrm{d}x\leq C\epsilon.
\end{equation}
By $(F_1)$, let $\lambda\triangleq\max\{p,q\}$, we can choose a constant $C_\epsilon>0$ such that
\begin{equation*}
  F(s)\leq C_\epsilon|s|^\lambda,\quad\text{for all} \,\,\, s\leq s_\epsilon.
\end{equation*}
Then, for $R>0$, we have
\begin{equation*}
  \begin{split}
    \int_{\{|x|>R\}\bigcap\{|u_n|\leq s_\epsilon\}}F(u_n)\mathrm{d}x
    & \leq C_\epsilon \int_{\{|x|>R\}\bigcap\{|u_n|\leq s_\epsilon\}}|u_n|^\lambda\mathrm{d}x \\
    & \leq C_\epsilon \int_{\{|x|>R\}\bigcap\{|u_n|\leq s_\epsilon\}}|u_n-u|^\lambda\mathrm{d}x+C_\epsilon \int_{\{|x|>R\}\bigcap\{|u_n|\leq s_\epsilon\}}|u|^\lambda\mathrm{d}x.
  \end{split}
\end{equation*}
Now using the compactness of embedding $X_\beta\hookrightarrow L^\lambda(\mathbb{R}^4)$, we can choose $R>0$ such that
\begin{equation}\label{Fun2}
   \int_{\{|x|>R\}\bigcap\{|u_n|\leq s_\epsilon\}}|F(u_n)|\mathrm{d}x\leq C\epsilon.
\end{equation}
Similarly, one has
\begin{equation}\label{Fun3}
   \int_{|u_n|\geq s_\epsilon}|F(u)|\mathrm{d}x\leq C\epsilon \quad \text{and} \quad \int_{\{|x|>R\}\bigcap\{|u_n|\leq s_\epsilon\}}|F(u)|\mathrm{d}x\leq C\epsilon.
\end{equation}
Noticing that
\begin{equation}\label{Fun4}
  \int_{\{|x|<R\}\bigcap\{|u_n|\leq s_\epsilon\}}|F(u_n)-F(u)|\leq\int_{|x|<R}|F(u_n)-F(u)|<\epsilon.
\end{equation}
Combining (\ref{Fun1}), (\ref{Fun2}), (\ref{Fun3}) and (\ref{Fun4}), we have the fact that $F(u_n)\rightarrow F(u)$ in $L^1(\mathbb{R}^4)$.

\noindent Furthermore, we get $I_\beta'(u)=0$, namely, $u$ is a weak solution of (\ref{problem}).

If $u=0$,
since
\begin{equation*}
  \int^{}_{\mathbb{R}^4}F(u_n)\mathrm{d}x\rightarrow\int^{}_{\mathbb{R}^4}F(0)\mathrm{d}x=0 \quad \text{and} \quad I_{\beta}(u_n)\rightarrow m,
\end{equation*}
then
\begin{equation*}
  \|u_n\|_{\beta}\rightarrow(2m)^{\frac{1}{2}} < \left(\frac{\alpha_{\beta}}{\alpha}\right)^{\frac{1-\beta}{2}}.
\end{equation*}
By Lemma \ref{convergence1}, we deduce that $u_n\rightarrow0$ strongly in $X_\beta$ which is contradiction with the fact that $m>0$.

Thus, $u\neq0$.

Clearly, up to a subsequence, we have
\begin{equation*}
  \liminf_{n\rightarrow+\infty}\|u_n\|_\beta \geq\|u\|_\beta >0.
\end{equation*}
Set
\begin{equation*}
  v_n=\frac{u_n}{\|u_n\|_\beta},\, v=\frac{u}{\lim\limits_{n\rightarrow+\infty}\|u_n\|_\beta}.
\end{equation*}
Then $\|v\|_\beta<1$, $v_n\rightharpoonup v$ weakly in $X_\beta$, $\|v_n\|_\beta=1$.

\noindent We have
\begin{equation*}
  \begin{split}
     \lim_{n\rightarrow+\infty}\|u_n\|^2_{\beta} & =2\left(m+\int^{}_{\mathbb{R}^4}F(u)\mathrm{d}x\right) \\
     & =\frac{2\left(m+\int^{}_{\mathbb{R}^4}F(u)\mathrm{d}x\right)\left(m-I_{\beta}(u)\right)}{\left(m-I_{\beta}(u)\right)} \\
     & =2\left(m-I_{\beta}(u)\right) \frac{m+\int^{}_{\mathbb{R}^4}F(u)\mathrm{d}x}{\int^{}_{\mathbb{R}^4}F(u)\mathrm{d}x-\frac{1}{2}\|u\|^2_\beta} \\
     & =2\left(m-I_{\beta}(u)\right)\frac{1}{1-\|v\|^2_\beta} \\
     & \leq\frac{2m}{1-\|v\|^2_\beta}
       <\frac{\left( \frac{\alpha_\beta}{\alpha}\right)^{1-\beta}}{1-\|v\|^2_\beta}.
  \end{split}
\end{equation*}
It follows that there exist $r_1>1$ and $0<\gamma_1<1$ such that, up to a subsequence,
\begin{equation*}
  \alpha r_1\|u_n\|^{\frac{2}{1-\beta}}_\beta
  <\frac{\gamma_1\alpha_\beta}{\left(1-\|v\|^2_\beta\right)^{\frac{1}{1-\beta}}}=\gamma_1P_\beta(v)\alpha_\beta,\quad\text{for all} \, \,n,
\end{equation*}
where $P_\beta$ is defined in (\ref{therom2.1(2)}).

\noindent By Theorem \ref{conclu2}, we get
\begin{equation*}
  \begin{split}
     \sup\limits_n\int^{}_{\mathbb{R}^4}\left(e^{\alpha r_1|u_n|^{\frac{2}{1-\beta}}}-1\right)\mathrm{d}x
     & =\sup\limits_n\int^{}_{\mathbb{R}^4}\left(e^{\alpha r_1\|u_n\|^{\frac{2}{1-\beta}} |v_n|^{\frac{2}{1-\beta}}}-1\right)\mathrm{d}x \\
     & \leq \sup\limits_n\int^{}_{\mathbb{R}^4}\left(e^{\gamma_1P_\beta(v)\alpha_\beta|v_n|^{\frac{2}{1-\beta}}}-1\right)\mathrm{d}x
       <+\infty.
  \end{split}
\end{equation*}
Continuing exactly as in the proof of Lemma \ref{convergence1}, we can show that $u_n\rightarrow u$ strongly in $X_\beta$.
\end{proof}

\vspace{1em}
\noindent\textbf{Proof of Theorem \ref{conclu3}}
In view of Lemma \ref{mp}$-$Lemma \ref{end}, we can apply the Mountain-Pass Theorem to deduce that $I_\beta$ admits a critical point $u$ such that $I_\beta(u)=m>0$. The nonnegativity of $u$ is immediate.
From the proof of Lemma \ref{end}, we can infer that $u$ is nontrivial.

\subsection*{Acknowledgments}
The authors have been supported by National Natural Science Foundation of China 11971392 and Natural Science Foundation of Chongqing, China cstc2021ycjh-bgzxm0115.

\end{document}